\documentclass{aims}

\pdfoutput=1

\usepackage{amsmath}
  \usepackage{paralist}
  \usepackage{graphics} 
  \usepackage{epsfig} 
\usepackage{graphicx}  \usepackage{epstopdf}
 \usepackage[colorlinks=true]{hyperref}
\hypersetup{urlcolor=blue, citecolor=red}

\usepackage{tikz}
\usepackage{float}

\renewcommand\dim{{\rm dim}}
\renewcommand\div{{\rm div}}
\newcommand{\deltasw}{\Delta_{sw}\, }
\newcommand{\nablasw}{\nabla_{\! sw}\, }
\newcommand{\divsw}{\div_{\! sw}\, }

\def\currentvolume{X}
 \def\currentissue{X}
  \def\currentyear{200X}
   \def\currentmonth{XX}
    \def\ppages{X--XX}
     \def\DOI{10.3934/xx.xx.xx.xx}

\theoremstyle{definition}

\newtheorem{remark}{Remark}

\title[Simulation of a dispersive shallow water system] 
      {A combined finite volume - finite element scheme for a dispersive shallow water system}

\author[N. A\"issiouene,  M-O. Bristeau, E. Godlewski and J. Sainte-Marie ]{N. A\"issiouene, M-O. Bristeau, E. Godlewski and J. Sainte-Marie}


 \email{Nora.Aissiouene@inria.fr}
 \email{Marie-Odile.Bristeau@inria.fr}
 \email{Edwige.Godlewski@upmc.fr}
 \email{Jacques.Sainte-Marie@inria.fr}

\thanks{The first author is supported by NSF grant xx-xxxx}

\thanks{$^*$ Corresponding author: xxxx}

\begin{document}
\maketitle

\medskip
{\footnotesize
\centerline{\small{Inria, EPC ANGE, Rocquencourt- B.P. 105, F78153 Le Chesnay cedex, France}}
\centerline{\small {CEREMA, EPC ANGE, 134 rue de Beauvais, F-60280 Margny-Les-Compiegne , France}}
\centerline{\small {Sorbonne Universit\'es, UPMC Univ Paris 06, UMR
    7598, Laboratoire Jacques-Louis Lions, F-75005, Paris, France}}
\centerline{\small {CNRS, UMR 7598, Laboratoire Jacques-Louis Lions, F-75005, Paris, France}}
} 

\medskip


\bigskip

 \centerline{(Communicated by the associate editor name)}

\begin{abstract}

We propose a variational framework  for the resolution of a
 non-hydrostatic Saint-Venant type model with bottom topography. This model is a shallow
 water type approximation of the incompressible Euler system with free
 surface and slightly differs from the Green-Nagdhi model,
 see~\cite{JSM_nhyd} for more details about the model derivation.

The numerical approximation relies on a prediction-correction type scheme initially introduced by Chorin-Temam \cite{chorin} to treat the incompressibility in the Navier-Stokes equations. The hyperbolic part of
the system is approximated using a kinetic finite volume solver
and the correction step implies to solve a mixed problem where the velocity and the pressure are defined in compatible finite element spaces.

The resolution of the incompressibility constraint leads to an elliptic problem involving the non-hydrostatic part of the pressure. This step uses a variational formulation of a shallow water version of the incompressibility condition.

Several numerical experiments are performed to confirm the relevance of our approach.


\end{abstract}

\section{Introduction}

Starting from the incompressible Euler or Navier-Stokes system, the hydrostatic assumption
consists in neglecting the vertical acceleration of the fluid. More
precisely~-- and with obvious notations~-- the momentum along the
vertical axis of the Euler equation
$$\frac{\partial w}{\partial t} + u\frac{\partial w}{\partial x} +
w\frac{\partial w}{\partial z} + \frac{\partial p}{\partial z} = -g,$$
reduces in the hydrostatic context to
\begin{equation}
\frac{\partial p}{\partial z} = -g.
\label{eq:hydro}
\end{equation}
Such an assumption produces important consequences over the structure and
complexity of the model. Indeed, Eq.~\eqref{eq:hydro} implies that the
pressure $p$ is no longer the Lagrange multiplier of the
incompressibility constraint and $p$ can be expressed, for free
surface flows, as a
function of the water depth of the fluid. Therefore, the hydrostatic
assumption implies that the resulting model, even though it describes an
incompressible fluid, has common features with models arising
in compressible fluid mechanics.\\

In geophysical problems, the hydrostatic assumption coupled with a shallow
water type description of the flow is often used. Unfortunately, these models do not represent phenomena
containing dispersive effects for which the non-hydrostatic contribution
cannot be neglected. And more complex models have to be considered to
take into account this kind of phenomena, together with numerical methods able
to discretize the high order derivative terms coming from the dispersive effects. Many shallow water type dispersive models
have been proposed such as KdV, Boussinesq, Green-Naghdi,
see~\cite{green,camassa1,bbm,nwogu,peregrine,lannes_marche,duran_marche,bonneton_lannes,lannes,lannes1,JSM_DCDS}. The modeling of the non-hydrostatic effects for shallow water flows
does not raise insuperable
difficulties but their discretization is more tricky. Numerical techniques for the approximaion of these models have been
recently proposed~\cite{chazel2,JSM_CF,lemet-gavri}.

The model studied in the present paper has been derived and studied
in~\cite{JSM_nhyd}. Its numerical approximation based on a
projection-correction strategy~\cite{chorin} is described
in~\cite{JSM_nhyd_num}. In~\cite{JSM_nhyd_num}, the discretization of
the elliptic part arising from the non-hydrostatic terms is carried
out in a finite difference framework. It is worth noticing that the
numerical scheme given in~\cite{JSM_nhyd_num} is endowed with
robustness and stability properties such as positivity,
well-balancing, discrete entropy and wet/dry interfaces treatment.

The main contents of this paper
is the derivation and validation of the correction step in a
variational framework.
Since the derivation in a 2d context of the model proposed in~\cite{JSM_nhyd} does not raise
difficulty, the results depicted in this paper pave the
way for a discretization of the 2d model over an unstructured mesh, and we will often maintain general notations as far as possible.

Notice that the non-hydrostatic model we consider slightly differs from the
well-known Green-Naghdi model~\cite{green} but the numerical
approximation proposed in this paper can also be used for the
numerical approximation of the Green-Naghdi system.

Let $\Omega \subset \mathbb{R} $, be a 1d domain (an interval) and $\Gamma=\Gamma_{in}\cup \Gamma_{out} $ its boundary (see figure \ref{figdessin}). The non-hydrostatic model derived in ~\cite{JSM_nhyd,JSM_nhyd_num} reads
\begin{eqnarray}
\frac{\partial H}{\partial t} + \frac{\partial H\overline{u}}{\partial x} &=&0, \label{nh1}\\
\frac{\partial H\overline{u}}{\partial t} + \frac{\partial}{\partial x}(H\overline{u}^2 +g\frac{H^2}{2} + H\overline{p}_{nh}) &=& -(gH +2\overline{p}_{nh} )\frac{\partial z_b}{\partial x},\label{nh2}\\
\frac{\partial H\overline{w}}{\partial t} + \frac{\partial H\overline{w}\overline{u}}{\partial x} &=& 2\overline{p}_{nh},\label{nh3}\\
\frac{\partial H\bar{u}}{\partial x}-\bar{u}\frac{\partial (H+2 z_b)}{\partial x} +2\bar{w}&=& 0, \label{contrainte}
\end{eqnarray}
where $H$ is the water depth, $z_b$ the topography and $p_{nh}$ the
non-hydrostatic part of the pressure. The variables denoted with a bar recall that this model is obtained performing an average along the water depth of the incompressible Euler system with free surface. The velocity field is denoted
$\bold{\bar{u}}=(\bar{u},\bar{w})^t$ with $\bar{u}$ (resp.  $\bar{w}$)
the horizontal (resp. vertical) component.

We denote $\eta=H+z_b$ the free surface of the fluid. In addition, we give the following notations
\begin{eqnarray}
\bold{n} =\left(
\begin{array}{ccc}
n \\
0
\end{array}
\right),\label{normal}
\end{eqnarray}
with $n$ the unit outward normal vector at $ \Gamma$ (in 1d, $n=\pm$ 1), $ \bold{n}$ represents the unit outward normal vector of the domain covered by the fluid, namely $\Omega \times [z_b,\eta] $.
We also consider the gradient operator
\begin{eqnarray}
\nabla_0=\left(
\begin{array}{ccc}
\frac{\partial}{\partial x} \\
0
\end{array}
\right).\label{grad}
\end{eqnarray}

\begin{figure}[htbp]
\begin{center}
\caption{Notations and domain definition.}
\label{figdessin}
\end{center}
\end{figure}

The smooth solutions of the system~\eqref{nh1}-\eqref{contrainte} satisfy moreover an energy conservation law
\begin{eqnarray}
\frac{\partial \bar{E} }{\partial t}+\frac{\partial }{\partial x}\left(\bar{u}\left( \bar{E}+\frac{g}{2}H^2+H \bar{p}_{nh}\right)\right) =0,
\end{eqnarray}
with
\begin{eqnarray}
 \bar{E} = 	\frac{H\left( \bar{u}^2+\bar{w}^2\right)}{2}
+\frac{gH\left(\eta + z_b \right)}{2}.
\end{eqnarray}

Note that \eqref{contrainte}  represents a shallow water version of the divergence free constraint, for which the non hydrostatic pressure $\bar{p}_{nh}$ plays the role of a Lagrange multiplier. Notice that considering $\bar{p}_{nh}=0$ and neglecting~\eqref{nh3}, the system~\eqref{nh1}-\eqref{nh2},\eqref{contrainte} reduces to the
classical Saint-Venant system.


The paper is organized as follows. First we give a rewriting of the
model and we present the prediction-correction method, the main part
being the variationnal formulation of the correction part. Then in
Section~\ref{sec:num_approx}, we detail the numerical
approximation. Finally, in Section~\ref{sec:valid}, numerical simulations
validating the proposed discretization techniques are presented.

%

\section{The projection scheme for the non-hydrostatic model}

Projection methods have been introduced by A. Chorin and R. Temam \cite{chorin_proj} in order to compute the pressure for
incompressible Navier-Stokes equations. These methods, based on a time
splitting scheme, have been widely studied and applied
to treat the incompressibility constraint (see
\cite{GuermondS04,Shen92onpressure,Shen_pseudo-compressibilitymethods}). We develop below an analogue of this method for shallow water flow. In order to describe the fractional time step method we use, we propose a rewritting of the model
\eqref{nh1}-\eqref{contrainte}.

\subsection{A rewritting}

Let us introduce the two operators $\nablasw$ and $\divsw$ defined by
\begin{eqnarray}
\nablasw f &=&\left(
\begin{matrix}
H\frac{\partial f}{\partial x}
+f\frac{\partial( H + 2z_b)}{\partial x} \\
-2f
\end{matrix}
\right)\label{eq:gradsw},\\
\divsw (\bold{v})&=&\frac{\partial H v_1}{\partial x}-v_1\frac{\partial (H+2z_b)}{\partial x}+2v_2\label{eq:divsw},
\end{eqnarray}
with $\bold{v}=(v_1,v_2)^t$. We assume for a while that $f$ and $\bold{v}$ are smooth enough. The shallow water form of the divergence
operator $\divsw $ (resp. of the gradient operator $\nablasw$)
corresponds to a depth averaged version of the divergence
(resp. gradient) appearing in the incompressible Euler and
Navier-Stokes equations.
Notice that the two operators $\nablasw$, $\divsw$ defined
by~\eqref{eq:gradsw}-\eqref{eq:divsw} are $H$ and $z_b$
dependent and we assume that $H$ and $z_b$ are sufficiently smooth
functions. One can check that these operators verify the fundamental duality relation
\begin{eqnarray}
\int_\Omega \divsw (\bold{v})f \,dx = -\int_\Omega \nablasw f\cdot \bold{v} \,dx + \left[ Hv_1\,f\right]_\Gamma.  \label{duality}
\end{eqnarray}
These definitions allow to rewrite the model \eqref{nh1}-\eqref{contrainte} as
\begin{eqnarray}
\frac{\partial H}{\partial t} +\frac{\partial H\bar{u}}{\partial x}&=&0 \label{system1} ,\\
\frac{\partial H\bold{\bar{u}} }{\partial t}
+\frac{\partial}{\partial x}\left( \overline{u}
  H\overline{\bf u}\right)
+ \nabla_0\left( \frac{g}{2}H^2 \right)
+\nablasw \bar{p}_{nh} &=&-gH\nabla_0 z_b \label{system2},\\
\divsw (\bold{\bar{u}})&=&0 ,\label{system3}
\end{eqnarray}
with $\nabla_0$ defined by~\eqref{grad}.

The system \eqref{system1}-\eqref{system3} can be written in the
compact form
\begin{eqnarray}
\frac{\partial X}{\partial t}+ \frac{\partial }{\partial x}F(X)+R_{nh}
&=&S(X),\label{F1} \\
\divsw (\bold{\bar{u}})&=&0,\label{F2}
\end{eqnarray}
where we denote
\begin{eqnarray}
 X
=\left( \begin{array}{ccc}
H\\
H\bar{u} \\
H\bar{w}\\
\end{array}\right),
\quad\quad
F(X)=
\left(\begin{array}{ccc}
H\bar{u} \\
H\bar{u}^2+\frac{g}{2}H^2 \\
H\bar{u}\overline{w}\\
\end{array}\right),
\end{eqnarray}
and
\begin{eqnarray}
R_{nh}=\left(\begin{array}{ccc}
0 \\
\nablasw \bar{p}_{nh}
\end{array}\right),
\quad\quad
S(X)=\left(\begin{array}{ccc}
0 \\
-gH\nabla_0 z_b
\end{array}\right).
\end{eqnarray}

Let be given time steps $\Delta t^n$ and note $t^n=\sum_{k\leq n} \Delta
t^k$. As detailed in~\cite{JSM_nhyd_num},
the projection scheme for system \eqref{F1}-\eqref{F2} consists in the following time splitting
\begin{eqnarray}
	    X^{n+1/2} &=& X^n- \Delta t^n \,\frac{\partial}{\partial x} F(X^n)+\Delta t^n\, S(X^n), \label{frac1}\\
	    X^{n+1}   &=& X^{n+1/2}- \Delta t^n\, R_{nh}^{n+1}    ,       \label{frac2} \\
	    \divsw {\bold{\bar{u}}^{n+1}}&=&0,\label{frac3}
	\end{eqnarray}
with $\bold{\bar{u}}^{n+1}=\left(
  \frac{(H\bar{u})^{n+1}}{H^{n+1}},\frac{(H\bar{w})^{n+1}}{H^{n+1}}\right)^t$.

\noindent The first two equations of~\eqref{frac1} consist in the
classical Saint-Venant system with topography and the third equation
is an advection equation for the quantity $H\bar{w}$. Equations~\eqref{frac2}-\eqref{frac3} describe the correction
step allowing to determine the non hydrostatic part of the pressure
$p_{nh}^{n+1}$ and hence giving the corrected state $X^{n+1}$. The numerical
resolution of~\eqref{frac1}~-- especially the first two equations~--
has received an extensive
coverage and efficient and robust numerical techniques exist, mainly
based on finite volume approach,
see~\cite{bouchut_book,JSM_entro}. The derivation of a robust and
efficient numerical technique for the resolution of the correction
step \eqref{frac2}-\eqref{frac3} is the key point. A strategy based on
a finite difference approach has been proposed, studied and validated
in~\cite{JSM_nhyd_num}. Unfortunately, the finite difference framework
does not allow to tackle situations with unstructured meshes in 2 or 3
dimensions. It is the key point of this paper to propose a variational
formulation of the correction step coupled with a finite volume
discretization of the prediction step.\\

The time discretization in the numerical scheme described above
corresponds to a fractional time step strategy with a first order Euler
scheme, explicit for the hyperbolic part and implicit for the elliptic
part. For hyperbolic conservation laws, the second-order accuracy in time is usually recovered by the Heun
method \cite{heun,bouchut_book} that is a slight modification of the second order
Runge-Kutta method. More precisely, for a dynamical system written
under the form
\begin{equation}
\frac{\partial y}{\partial t} = F(y),
\label{eq:sysd}
\end{equation}
the Heun scheme consists in defining $y^{n+1}$ by
\begin{equation}
y^{n+1}= \frac{y^n+\tilde{y}^{n+2}}{2},
\label{eq:heun}
\end{equation}
with
\begin{eqnarray}
& & \tilde{y}^{n+1} = y^n + \Delta t^n F(y^n),\nonumber\\
& & \tilde{y}^{n+2} = \tilde{y}^{n+1} + \Delta t^n F(\tilde{y}^{n+1}).\nonumber
\end{eqnarray}
The model we have to discretize has the general form
\begin{align}
& \frac{\partial X}{\partial t} = F(X,p),\label{eq:eq_sys1}\\
& B(X) = 0,\nonumber
\end{align}
and we propose the numerical scheme
\begin{equation}
X^{n+1}  = \frac{X^n+\tilde{X}^{n+2}}{2},
\label{eq:heun_mod}
\end{equation}
with the two steps defined by
\begin{eqnarray}
& & \tilde{X}^{n+1} = X(t^n) + \Delta t^n F(X^n,\tilde{p}^{n+1}),\label{eq:2nd_order11}\\
& &  B(\tilde{X}^{n+1}) = 0 \label{eq:div_2nd_order},
\end{eqnarray}
and
\begin{eqnarray}
& & \tilde{X}^{n+2} = \tilde{X}^{n+1} + \Delta t^n
F(\tilde{X}^{n+1},\tilde{p}^{n+2}),\label{eq:2nd_order2}\\
& &  B(\tilde{X}^{n+2}) = 0 \label{eq:div_2nd_order2},
\end{eqnarray}
where $\tilde{p}^{n+1}$ and $\tilde{p}^{n+2}$ are the solutions of the elliptic equations
derived from the divergence free constraints~\eqref{eq:div_2nd_order}
and~\eqref{eq:div_2nd_order2} respectively.
Notice that \eqref{eq:2nd_order11} is a compact form of the fractional scheme \eqref{frac1}-\eqref{frac2} where the intermediate step no more appears.

\noindent It is easy to prove the scheme~\eqref{eq:heun_mod} is second order
accurate. Indeed assuming $F$ and $B$ are enough smooth, we have
$$X(t^n+\Delta t^n) = X(t^n) + \Delta t^n \dot{X}(t^n) + \frac{(\Delta t^n)^2}{2} \ddot{X}(t^n) + {\mathcal{O}} \bigl((\Delta t^n)^3\bigr),$$
or equivalently using~\eqref{eq:eq_sys1}
\begin{multline}
X(t^n+\Delta t^n) = X(t^n) + \Delta t^n F(X(t^n),p(t^n)) \\
+ \frac{(\Delta
  t^n)^2}{2} \left(F(X(t^n),p(t^n))\frac{\partial F(X(t^n),p(t^n))}{\partial X}
  + \frac{\partial F(X(t^n),p(t^n))}{\partial p}\frac{\partial p(t^n)}{\partial t}
\right) + {\mathcal{O}} \bigl((\Delta t^n)^3\bigr). \label{eq:dev1}
\end{multline}
Now from~\eqref{eq:2nd_order11} and~\eqref{eq:2nd_order2}, we get
\begin{eqnarray}
\tilde{X}^{n+2} & = & X^n + \Delta t^n F(X^n,\tilde{p}^{n+1}) + \Delta t^n
F(\tilde{X}^{n+1},\tilde{p}^{n+2})\nonumber\\
& = & X^n + 2\Delta t^n F(X^n,\tilde{p}^{n+1}) \nonumber\\
& & + (\Delta t^n)^2 \left(
  F(X^n,\tilde{p}^{n+1})\frac{\partial
    F(X^n,\tilde{p}^{n+1})}{\partial X} + \frac{\partial
    F(X^n,\tilde{p}^{n+1})}{\partial p}\frac{\partial
    \tilde{p}^{n+1}}{\partial t}\right) \nonumber \\
 &&  + {\mathcal{O}} \bigl((\Delta t^n)^3\bigr).
\label{eq:dev2}
\end{eqnarray}
Using~\eqref{eq:heun_mod}, we see that relations~\eqref{eq:dev1} and~\eqref{eq:dev2} are equivalent up to third order terms.


\subsection{The correction step}

In this part, we consider we have at our disposal a space
discretization of Eq.~\eqref{frac1} solving the hydrostatic part of
the model and we focus on the correction step~\eqref{frac2}-\eqref{frac3}.

\subsubsection{Variational formulation}
The correction step~\eqref{frac2}-\eqref{frac3} writes,
			\begin{eqnarray}
				H^{n+1}&=&H^{n+1/2},  \label{proj1}\\
			(H\bold{u})^{n+1}+ \Delta t^n \nablasw {p}_{nh}^{n+1}&=&(H\bold{u})^{n+1/2},  \label{proj2}  \\
		  \divsw (\bold{u}^{n+1})&=&0.		\label{proj3}
			\end{eqnarray}
For the sake of clarity, in the following we will drop the notation
with a bar and we denote $p$ instead of $\bar{p}_{nh}$. Likewise we drop
the superscript $^{n+1}$ for the corrected states.

Equations~\eqref{proj2}-\eqref{proj3} is a mixed problem in
velocity/pressure, its approximation leads to the variational mixed problem

Find ${p} \in Q$, $\bold{u} \in V_0$
with
\begin{eqnarray}
Q=\{ q  \in L^2(\Omega) |\nablasw q  \in (L^2(\Omega))^2 \},
\end{eqnarray}

and
\begin{eqnarray}
V_0=\{ {\bf v}=(v_1,v_2) \in (L^2(\Omega))^2 |\divsw  (\bold{v} ) \in L^2(\Omega), {v_1}_{|\Gamma}=0\},
\end{eqnarray}

such that
\begin{eqnarray}
\int_\Omega \left(H\bold{u} + \Delta t^n \nablasw p
 \right) \cdot \bold{v} \, dx &=&\int_\Omega (H\bold{u})^{n+1/2} \cdot \bold{v} \ dx,
\quad \forall \, \bold{v} \in V_0,\label{varia_1}\\
\int_\Omega \divsw (\bold{u}) q \, dx &=& 0,
\quad \forall q \in Q. \label{varia_2}
\end{eqnarray}

Introducing the bilinear forms
\begin{eqnarray*}
 a(\bold{u},\bold{v})&=&\int_{\Omega} H\bold{u} \cdot \bold{v} \, dx, \quad \forall \, \bold{u},\bold{v} \in V_0,\\
 b(\bold{v},q)&=&-\int_{\Omega} \divsw (\bold{v}) q \, dx,  \quad \forall \bold{v} \in V_0,\forall \, q \in Q,
\end{eqnarray*}
the problem \eqref{varia_1}-\eqref{varia_2} becomes:\\
Find $p \in Q $ and $\bold{u} \in V_0 $ such that
\begin{eqnarray}
\frac{1}{\Delta t^n}a(\bold{u},\bold{v})+b(\bold{v},p)
&=&\frac{1}{\Delta t^n}a(\bold{u}^{n+1/2},\bold{v}), \quad \forall
\bold{v}  \in V_0, \label{mixed_pb1} \\
b(\bold{u},q)&=&0 \label{mixed_pb2},\quad \forall q \in Q.
\end{eqnarray}

\subsubsection{The pressure equation}

Formally, we take $\bold{v}$ in the form $\bold{v}=\frac{\nablasw q}{H}$ with $  q \in Q_{0}$ and
$$
Q_0=\{ q  \in Q \, q|_{\Gamma} = 0 \}. $$

\noindent
Then, with  \eqref{duality} and \eqref{mixed_pb2}, we get
\begin{eqnarray*}
a(\bold{u,\bold{v}})&=&\int_{\Omega} u\cdot \nablasw q dx \\
&=&-\int_{\Omega} \divsw(\bold{u})\cdot q dx  \\
&=&0.
\end{eqnarray*}
If we introduce the shallow water version of the Laplacian operator $\deltasw$
\begin{eqnarray*}
\deltasw p=\divsw \left(\frac{1}{H}\nablasw p\right)\label{lap},
\end{eqnarray*}
it is natural to consider the new variational formulation:

Find $\bold{u} \in V, p \in Q_{0,sw}$ such that

\begin{eqnarray}
\frac{1}{\Delta t^n}a(\bold{u},\bold{v})
+(\bold{v},\nablasw p)
&=&\frac{1}{\Delta t^n}a(\bold{u}^{n+1/2},\bold{v}) \quad \forall
\bold{v} \in V, \label{mixLap1} \\
 \left(\deltasw p,q \right)
&=&
\frac{1}{\Delta t^n} \left( \divsw (\bold{u}^{n+1/2}),q\right), \quad \forall q \in Q_{0,sw},\label{mixLap2}
\end{eqnarray}
with  $Q_{0,sw}=\{ q\in Q_0 | \divsw\left( \frac{\nablasw q}{H}\right) \in L^2(\Omega) \} $.
\medskip

From \eqref{mixLap2}, we deduce
\begin{eqnarray}
\deltasw (p)=\frac{1}{\Delta t^n} \divsw (\bold{{u}}^{n+1/2})\label{ellipLap}.
\end{eqnarray}

\noindent The resolution of the equation \eqref{ellipLap} allows to update the velocity at the correction step \eqref{proj2}.

\noindent Notice that the equation \eqref{ellipLap} is equivalent to apply the operator $\divsw $ to the equation \eqref{frac2} and to use the shallow water free divergence \eqref{frac3} to eliminate $\bold{u}$.

\begin{remark}
Taking the functional spaces $ V$ and $ Q_{0,sw}$, the problem \eqref{mixLap1}-\eqref{mixLap2} is not equivalent to the problem \eqref{mixed_pb1}-\eqref{mixed_pb2}.
\end{remark}

\begin{remark}
Notice that Eq.~\eqref{ellipLap} has the form of a Sturm-Liouville
type equation.
\end{remark}

\subsubsection{Inf-sup condition}

To ensure that the saddle problem \eqref{mixed_pb1}-\eqref{mixed_pb2} is well posed, the Babuska-Brezzi \cite{Brezzi1974,pironneau1988methodes} condition
\begin{eqnarray*}
\exists \gamma >0, \quad \gamma < \inf_{q \in Q}\sup_{\bold{v} \in V_0}
\frac{b(\bold{v},q)}{\vert\vert \bold{v}\vert\vert_{V_0}\vert\vert q\vert\vert_{Q}},
\end{eqnarray*}
has to be satisfied. Denoting by $\mathcal{B}$ the weak operator defined by $\forall v \in V_0 $, $ \mathcal{B}v = b(v,q), \forall q \in Q $, we have

\begin{eqnarray*}
\ker \mathcal{B}^t&=& \{ q \in Q | \mathcal{B}^t q=0\} \\
&=& \{ q\in Q | \int_\Omega \nablasw q \cdot \bold{v} \, dx=0 \; \forall \bold{v} \in V_0 \}.
\end{eqnarray*}
We assume $H\geq H_0 >0$. Because of the positivity of $H$, it is obvious that the bilinear form $a$ is coercive.
Choosing $\bold{v} \in (L^2(\Omega))^2 $ and $ q $ such as $ \nablasw q
\in (L^2(\Omega))^2 $, it follows that $ \nablasw q=0 $, then $q=0$ and $ \ker \mathcal{B}^t ={0}$. Indeed, in contrast with Navier-Stokes equations for which the pressure is defined up to an additive constant, the non hydrostatic of the shallow water equations is fully defined.
 Therefore, the mixed problem \eqref{mixed_pb1}-\eqref{mixed_pb2} satisfies the inf-sup condition and admits a unique solution.

\subsection{Boundary conditions}
\label{secBD}
In this section, we still consider  that the hydrostatic part is provided and we study the compatibility of the boundary conditions between the hydrostatic part and the projection part. Therefore, the compatibility between the pressure and velocity at boundary needs to be studied.
To this aim, we first provide the conditions required to impose Dirichlet or Neumann pressure at boundary, and then, we couple these conditions with the hydrostatic part.

We consider a more general case taking the space $V$

$$ V=\{ \bold{v}=(v_1,v_2) \in (L^2(\Omega))^2 |\divsw  \bold{v}  \in L^2(\Omega) \},$$

\noindent and we introduce the bilinear form
$$c(\bold{v},p)=\int_{\Gamma} Hp \, v_1 \,n \, ds,  \quad \forall \bold{v} \in V, p \in Q, $$
with $n$ the unit outward normal vector defined by~\eqref{normal}. In one dimension, $c(\bold{v},p)= (Hpv_1)|_{\Gamma_{out}}-(Hpv_1)|_{\Gamma_{in}} $.\\

Therefore instead of \eqref{mixed_pb1}-\eqref{mixed_pb2}, we consider the problem:\\
Find $ \bold{u} \in V$, $p\in Q$ such that,
\begin{eqnarray}
\frac{1}{\Delta t^n}a(\bold{u},\bold{v})+b(\bold{v},p)
&=&
\frac{1}{\Delta t^n}a(\bold{u}^{n+1/2},\bold{v})+c(\bold{v},p),  \quad \forall \, \bold{v}  \in V,\label{mixed_bc1}\\
b(\bold{u},q)&=&0, \quad \forall q \in Q.\label{mixed_bc2}
\end{eqnarray}
Notice that $ \divsw (\bold{u})=\nabla_0\cdot(H\bold{u})+\bold{u}\cdot(\bold{n_s}+\bold{n_b})$ and $\nablasw p=H\nabla_0(p)-p(\bold{n_s}+\bold{n_b}) $ with $\nabla_0$ defined by \eqref{grad} and $ \bold{n_s}$ (resp. $ \bold{n_b}$) the (non-unit) normal vector at the surface (resp. at the bottom)
\begin{eqnarray*}
\bold{n_s}=
\left(\begin{array}{ccc}
-\frac{\partial \eta }{\partial x}\\
1
\end{array}\right),
\qquad
\bold{n_b}=
\left(\begin{array}{ccc}
-\frac{\partial z_b}{\partial x} \\
1
\end{array}\right).
\end{eqnarray*}
Moreover, we have
\begin{eqnarray*}
\int_{\Omega} \divsw (\bold{u})\;dx=\int_{\Gamma} H{u} \, n \, ds+\int_{\Omega} \bold{u}\cdot (\bold{n_s}+\bold{n_b}) \, dx.
\end{eqnarray*}
Hence, to satisfy the divergence free condition, the velocity $ \bold{u}$ should verify
\begin{eqnarray*}
\int_{\Gamma} H{u} \,n \, ds=-\int_{\Omega} \bold{u}\cdot (\bold{n_s}+\bold{n_b}) \, dx .
\end{eqnarray*}

\subsubsection*{Dirichlet condition for the pressure}

From the variational
formulation \eqref{mixed_pb1}-\eqref{mixed_pb2} of the projection
scheme, a natural boundary condition for the pressure is a Dirichlet
condition. At  $ \Gamma_i, (i=in,out)$, $p|_{\Gamma_i}=p_0 $ then
$c(\bold{v},p)=\int_{\Gamma_i} p_0v_1  n\, ds$ and we take $ \bold{v}
\in V$, $q\in Q_i$ with
\begin{eqnarray*}
Q_i&=&\{ q  \in L^2(\Omega) |\nablasw  q  \in (L^2(\Omega))^2 , q|_{\Gamma_{i}}=0\}.
\end{eqnarray*}

\subsubsection*{Neumann boundary for the pressure}

The Neumann boundary condition for the projection scheme is not natural and to enforce such a condition, the elliptic  problem \eqref{mixLap1}-\eqref{mixLap2} has to be considered.
Taking now $q \in Q_{sw}$, with $ Q_{sw}=\{ q \in Q| \divsw(\frac{\nablasw p}{H}) \in L^2(\Omega)^2 \}$,
the problem is rewritten
\begin{eqnarray*}
\frac{1}{\Delta t^n}a(\bold{u},\bold{v})+(\bold{v},\nablasw p)
&=&\frac{1}{\Delta t^n}a(\bold{u}^{n+1/2},\bold{v})
 \quad \forall
\bold{v} \in V,  \\
 \left(\deltasw p,q \right)
&=&
\frac{1}{\Delta t^n} \left( \divsw (\bold{u}^{n+1/2}),q\right)+\frac{1}{\Delta t^n}\tilde{c}(\bold{u},q), \quad \forall q \in Q_{sw},
\end{eqnarray*}

with $\tilde{c}$, the bilinear form
$$\tilde{c}(\bold{u},p)=
\int_{\Gamma}  \left( Hu
+ \Delta t^n \nabla_{sw}p|_1 \, \right) q \, n \,d \gamma
- \int_{\Gamma}  Hu^{n+1/2} \, q \, n \, d \gamma. $$
Many studies have been done to choose an appropriate variational formulation for this problem. In \cite{Guermond1996} J-L. Guermond explores the different variational formulations in order to enforce a Neumann pressure boundary condition, in \cite{Johnston} some equivalent formulations are given to switch between Neumann and Dirichlet boundary conditions.

Taking the normal component  at the boundary $ \Gamma_i $ of the momentum equation, it follows that
\begin{eqnarray*}
H\frac{\partial p}{\partial n}|_{\Gamma_i}
+p|_{\Gamma_i}(\frac{\partial H}{\partial n}|_{\Gamma_i})
=\frac{H}{\Delta t^n}
({u}|_{\Gamma_i}^{n+1/2}-{u}|_{\Gamma_i}) .
\end{eqnarray*}

We note $ \frac{\partial H}{\partial n}|_{\Gamma_i}=\beta_i,\, i=in,out$.
\begin{itemize}
\item 	In case $\beta_i =0$, a Neumann boundary condition for the pressure is deduced of a Dirichlet condition for $u$.
\begin{eqnarray}
\frac{\partial p}{\partial n}|_{\Gamma_i}&=&\frac{1}{\Delta t^n}({u}|_{\Gamma_i}^{n+1/2}-u|_{\Gamma_i}),\label{neum1}
\end{eqnarray}

\item In the other cases, it gives a mixed boundary condition

\begin{eqnarray}
\frac{\partial p }{\partial n}|_{\Gamma_i}+ \beta_i p  &=&\frac{1}{\Delta t^n}({u}|_{\Gamma_i}^{n+1/2} -u|_{\Gamma_i}).\label{neum2} 
\end{eqnarray}

\end{itemize}
Then, in the two cases, we have imposed a Dirichlet velocity
condition, that leads to take $\bold{v} \in V_i$ and $ q \in Q$, with for $i=in,out$
\begin{eqnarray}
V_i=\{ {\bf v}=(v_1,v_2) \in (L^2(\Omega))^2 |\divsw  (\bold{v} ) \in L^2(\Omega) , {v_1}_{|\Gamma_i}=0\}. \label{vi}
\end{eqnarray}
Let us now give the  coupling boundary conditions between the prediction step and the correction step. Indeed, in the projection part, boundary conditions need to be set in order to be consistent with the hydrostatic part.

Concerning the prediction step, we consider the well known Saint-Venant system and we assume that the Riemann invariant remains constant along the associated characteristic. This approach has been introduced in ~\cite{cdl_mob} and distinguishes fluvial and torrential boundaries depending on the Froude number $Fr=\frac{\left|{u}\right|}{c}$. Usual boundary conditions consist to impose a flux $\bold{q}_0$ at the inflow boundary and a water depth at the outflow boundary. It is also classical to let a free outflow boundary, setting a Neumann boundary condition for the water depth and for the velocity. In both cases, we give the boundary conditions that have to be set in the correction step.

\noindent
We consider the first situation in which we set a flux at the inflow  $\Gamma_{in}$ and a given depth at the outflow  $\Gamma_{out}$.
Assuming a fluvial flow, this case consists in solving a Riemann problem at the interface $\Gamma_{in} $ where the global flux is given by $ \bold{q}_0=(q_{01},q_{02})^t=(Hu^{n+1/2},Hw^{n+1/2})^t$. That gives the boundary values $H_0= H^{n+1/2}_0$, $u_0=\frac{q_{02}}{H^{n+1/2}_0} $ and $w_0=\frac{q_{02}}{H^{n+1/2}_0}$ from the hyperbolic part. This leads to obtain a Dirichlet condition for the pressure at the left boundary of the correction part.
\noindent \\
Moreover, if $H$ is given for the outflow, we preconize to give a mixed condition for the pressure that corresponds to the boundary condition \eqref{neum2}
\begin{eqnarray*}
\left.p \right|_{\Gamma_{in}} & = & 0,\\
\left. \frac{\partial p }{\partial n}\right|_{\Gamma_{out}} +\left.p
  \frac{\partial H}{\partial n}\right|_{\Gamma_{out}} & = & 0,
\end{eqnarray*}
that leads to take $\bold{u} \in V_{out} $, with the definition \eqref{vi} and
\begin{eqnarray*}
p \in Q_{in}&=&	\{q \in L^2(\Omega) | \nablasw q \in L^2(\Omega)    ,\quad q_{|\Gamma_{in}}=0 \}.\nonumber
\end{eqnarray*}

\noindent
We now consider the second situation in which we still impose a flux $\bold{q}_0 $ at the inflow and we set a free outflow boundary. In this case, we assume the two Riemann invariants are constant along the outgoing characteristics of the hyperbolic part (see \cite{cdl_mob}), therefore, we have a Neumann boundary condition for $H^{n+1/2}$ and ${u}^{n+1/2}$.
\begin{eqnarray*}
\left.\frac{\partial H}{\partial n}\right|_{\Gamma_{out}}=0,
\qquad
\left.\frac{\partial \bold{u}^{n+1/2}}{\partial n}\right|_{\Gamma_{out}}=0.
\end{eqnarray*}
Preserving these conditions at the correction step, it gives a Neumann boundary condition for the pressure of type \eqref{neum1}
\begin{eqnarray*}
\left.\frac{\partial p}{\partial n}\right|_{\Gamma_{out}}=0.
\end{eqnarray*}

\noindent
For an inflow given, the functional spaces will be defined by
$$\bold{u} \in V_{out}, p \in Q_{in}.
$$

\section{Numerical approximation}
\label{sec:num_approx}

\subsection{Discretization}
\label{disc}

This section is devoted to the numerical approximation and mainly for the correction step.
Let us be given a subdivision of $\Omega$ with $N$ vertices $x_1<x_2<...<x_N $ and we define the space step  $\Delta x_{i+1/2}=x_{i+1}-x_i $. We also note $\Delta x_i = x_{i+1/2}-x_{i-1/2}$ with $x_{i+1/2} = \frac{x_i+x_{i+1}}{2}$.

\subsubsection*{Prediction part}

For the prediction step~\eqref{frac1} i.e the hydrostatic part of the model, we use a finite volume scheme.
We introduce the finite volume cells $C_i$ centered at vertices $x_i$ such that  $\Omega=\cup_{i=1,N} C_i$.
Then, the approximate solution $X_i^n$ at time $t^n$
\begin{eqnarray*}
X_i^n \approx \frac{1}{\Delta x_i}\int_{C_i}X(x,t^n)dx,
\end{eqnarray*}
is solution of the numerical scheme
\begin{eqnarray*}
X_i^{n+1}=X_i^n-\sigma_i^n \left( \mathcal{F}^n_{i+1/2} - \mathcal{F}^n_{i-1/2} \right) + \sigma_i^n S_i^n,
\end{eqnarray*}
where $ \sigma_i^n = \frac{\Delta t^n}{\Delta x_i}$ and $\mathcal{F}$
(resp. $\mathcal{S}$) is a robust and efficient discretization of the
conservative flux $F(X)$ (resp. the source term $S(X)$). The time step
is determined through a classical CFL condition. Many numerical fluxes
and discretizations are available in the literature~\cite{bouchut_book,godlewski_book,leveque_book}, we choose
a kinetic based solver~\cite{JSM_entro} coupled with the hydrostatic reconstruction technique~\cite{bristeau1}.

\subsubsection*{Correction part}

Concerning the correction step \eqref{frac2}-\eqref{frac3}, we consider the discrete problem corresponding to the mixed problem \eqref{mixed_pb1}-\eqref{mixed_pb2}. We approach ($V_0$,$Q$) by the finite dimensional spaces ($V_{0h}$,$Q_h$)  and we note
$$N = \dim(V_{0h}), \qquad M = \dim(Q_h).$$

We also denote by $ \left(\varphi_{i}\right)_{i=1,N} $ and ${\left(\phi_{l}\right)_{l=1,{M}}}$ the basis functions of $ V_{0h}$ and $Q_h$ respectively. The finite dimensional spaces will be specified later on. We approximate $ (\bold{u},p) \in (V_0,Q)$ by $ (\bold{u}_h,p_h)\in(V_{0h},Q_h) $ such that
$$
\bold{u}_h(x) = \sum_{i=1}^{N}
\left(\begin{array}{ccc}
u_{i} \\
w_{i}
\end{array}\right)
\varphi_{i}(x),\qquad
p_h(x) = \sum_{l=1}^{M} p_{l} \, \phi_{l}(x).
$$


Therefore, we consider the discrete problem:\\ \\
Find $ \bold{u}_h \in V_{0h}$, $p_{h} \in Q_h$ such that
\begin{eqnarray}
\frac{1}{\Delta t^n}a(\bold{u}_h,\bold{v}_h)+b(\bold{v}_h,p_h)
&=&\frac{1}{\Delta t^n}a(\bold{u}_h^{n+1/2},\bold{v}_h), \quad \forall
\bold{v}_h  \in V_{0h} , \label{mixed_pb1_disc} \\
b(\bold{u}_h,q_h)&=&0 ,\label{mixed_pb2_disc}\quad \forall q_h \in Q_h,
\end{eqnarray}

Let us introduce the mass matrix $M_H$ given by
\begin{eqnarray*}
M_H = \left(\int_\Omega H\varphi_i \varphi_j dx \right)_\text{$1\leq i,j \leq N$},
\end{eqnarray*}
and the two matrices $B^t$, $B$
defined by
\begin{eqnarray*}
B^t= \left( \int_\Omega \nabla_{sw}{(\phi_{l})} \varphi_i dx \right)_{1\leq l \leq M,1\leq i \leq N}
,
B= -\left( \int_\Omega \div_{sw}(\varphi_j)\phi_{l} dx \right)_{1\leq l \leq M,1\leq j \leq N},
\end{eqnarray*}

and we denote
$$U = \begin{pmatrix}
u_1\\
\vdots\\
u_N\\
w_1\\
\vdots\\
w_N
\end{pmatrix},
\qquad
P = \begin{pmatrix}
p_1\\
\vdots\\
p_M
\end{pmatrix},
$$
Therefore, the problem \eqref{mixed_pb1_disc}-\eqref{mixed_pb2_disc} becomes
\begin{eqnarray*}
\left(
\begin{array}{ccc}
 \frac{1}{\Delta t^n} A_H && 	 \,  B^t  \\
 B && 0
\end{array}
\right)
\left(
\begin{array}{ccc}
 U \\
 P
\end{array}
\right)
=
\left(
\begin{array}{ccc}
 \frac{1}{\Delta t^n}A_HU^{n+1/2} \\
 0\\
\end{array}
\right),
\end{eqnarray*}
with
\begin{eqnarray*}
A_H=
\left(
\begin{array}{ccc}
M_H & 0   \\
0   & M_H
\end{array}
\right).
\end{eqnarray*}
Assuming that $M_H$ is invertible and eliminating the velocity $U$, we obtain the  following equation
\begin{eqnarray}
BA_H^{-1}B^tP=\frac{1}{\Delta t^n} B U^{n+1/2},\label{lapsw}
\end{eqnarray}
that is a discretization of the elliptic equation~\eqref{ellipLap} of Sturm-Liouville
type governing the pressure $p$.

We now take into consideration the boundary conditions in the more general problem \eqref{mixed_bc1}-\eqref{mixed_bc2}. The velocity $ \bold{u}$ is approximated by $\bold{u}_h \in  V_h$,  and the discrete problem is then written:

Find $ (\bold{u}_h,p_h) \in ( V _h,Q_h)$ such that
\begin{eqnarray*}
\frac{1}{\Delta t^n}a(\bold{u}_h,\bold{v}_h)+b(\bold{v}_h,p_h)
&=&\frac{1}{\Delta t^n}a(\bold{u}_h^{n+1/2},\bold{v}_h) +\frac{1}{\Delta t^n}c(\bold{u}_h,p_h)\,, \quad \forall
\bold{v}_h  \in V_{h} ,  \\
b(\bold{u}_h,q_h)&=&0 \quad \forall q_h \in Q_h.
\end{eqnarray*}

Considering the matrix $\Delta t^n C=\left( c\left(\varphi_i,\phi_l\right)\right)_{1\leq l \leq M,1\leq i \leq N}$ that contains the boundary terms, the equation \eqref{lapsw} becomes
\begin{eqnarray*}
BA_H^{-1}(B^t-C)P=\frac{1}{\Delta t^n} B U^{n+1/2},\label{lapsw1}
\end{eqnarray*}

This approach is suitable for the finite element approximation that is given
in the next section. However, it implies to inverse a mass matrix $
M_H$ that is not diagonal and depends on the water depth $H$. In
practice, we use the mass lumping technique introduced by Gresho  (\cite{Gresho_Chan_1988}) to avoid inverting the mass matrix in projection methods for Navier-Stokes incompressible system.

\subsection{Finite element $\mathbb{P}_1$/$\mathbb{P}_0$}

In this part, the problem is solved by the mixed finite element approximation  $\mathbb{P}_1/\mathbb{P}_0$ (see \cite{pironneau1988methodes}) on the domain $ \Omega = \cup_{l=1}^{M}K_l$ ( $M$ = $ N-1$ with $N$ the number of nodes), where the velocity is approximated by a continuous linear function and the pressure is approximated by a discontinuous piecewise constant function over each element

 $$ {\bold{u}_h \in V_h=\{ \bold{v}_h \in (\mathcal{C}^0(\Omega))^2 \mid \ \bold{v}_h|_{K_l} \in \mathbb{P}_1^2, \forall l=1,\dots,N-1 \} },$$ and
$$ {{p}_{h} \in Q_h =\{ q_h \mid \left. q_h \right|_{K_l} \in \mathbb{P}_0 \; , \; \forall l=1,\dots, M-1 \}}.$$
Using the discretization given in \ref{disc}, we denote by $K_{i+1/2}$
the finite element $[x_{i},x_{i+1}]$, then the pressure is constant on the finite element $K_{i+1/2}$.

For the sake of clarity, in this situation, let $\left(\phi_{j+1/2}\right)_{1\leq j \leq M}$ be the basis functions for the pressure $p_h$, and
$\left(\varphi_i\right)_{1 \leq i \leq N}$  the basis functions
for the velocity $\bold{u}_h$ such that
\begin{eqnarray*}
 \bold{u}_h(x) = \sum_{i=1}^{N}
\left(\begin{array}{ccc}
u_{i} \\
w_{i}
\end{array}\right)
\varphi_{i}(x)
,\quad
p_{h}(x) = \sum_{j=1}^{M} p_{j+1/2} \, \phi_{j+1/2}(x).
\end{eqnarray*}


We note $ \zeta = H+2z_b$ and assume $\zeta$ is approximated by a piecewise linear function $ \zeta_h$, namely $ \zeta_h(x)=\sum_{i=i}^N\zeta_i\varphi_i(x)$. We also note $ \frac{\partial \zeta_h}{\partial x}|_{i+1/2}=\frac{\zeta_{i+1}-\zeta_i}{\Delta x_{i+1/2}} = \chi_{i+1/2}$ the constant gradient of $ \zeta_h$ on the element $ K_{i+1/2}$.\\
 If we denote $ \underline{\varphi} = \left(\varphi,\varphi \right)^t$, then the  shallow water gradient operator is written

\begin{eqnarray*}
\int_\Omega \nabla_{sw}(p_h)\cdot \underline{\varphi_i}\,dx=
\left(
\begin{array}{ccc}
H_i(p_{i+1/2}-p_{i-1/2}) +
\frac{p_{i-1/2}}{2}\chi_{i-1/2}
+
\frac{p_{i+1/2}}{2}\chi_{i+1/2}\\
-(\Delta x_{i+1/2} p_{i+1/2}+\Delta x_{i-1/2} p_{i-1/2})
\end{array}
\right).
\end{eqnarray*}
\\
\noindent Similarly, the shallow water divergence operator writes
\begin{eqnarray}
\int_\Omega \divsw(\bold{u}_h)\phi_{j+1/2}\,dx
=
{u}_{j+1}-{u}_j
+\frac{{u}_j+{u}_{j+1}}{2}\left(\zeta_{j+1}-\zeta_j \right)+\Delta x_{j+1/2}({w}_j+{w}_{j+1}). \nonumber
\end{eqnarray}

In one dimension, this approach corresponds to a staggered-grid finite-difference method where the velocity is computed at the nodes and the pressure is computed at the middle nodes. The discretization we obtain corresponds exactly to the finite difference scheme given in \cite{JSM_nhyd_num}, and then, the properties established in \cite{JSM_nhyd_num} are conserved.

\subsection{Finite element $\mathbb{P}_1$-iso-$\mathbb{P}_2$/$\mathbb{P}_1$}

For the one dimensional $\mathbb{P}_1$-iso-$\mathbb{P}_2/\mathbb{P}_1$, we consider two meshes $\mathcal{K}_{h}$ (the same as before) and $ \mathcal{K}_{2h}$  with $K_{h,i+1/2}=[x_i,x_{i+1}]$ and  $K_{2h,j}=[x_{2j-1},x_{2j+1}]$  the finite elements defined on the respective meshes $ \mathcal{K}_h$ and $ \mathcal{K}_{2h}$ such that $ \mathcal{K}_h=\cup_{i=1}^{N-1} K_{h,i+1/2}$ and $ \mathcal{K}_{2h}=\cup_{j=1}^{M-1} K_{2h,j}$ with $ N$ the total number of vertices of $\mathcal{K}_{h} $ and $ M= (N-1)/2$ (assuming $N$ odd), the number of vertices of $\mathcal{K}_{2h} $.
Therefore, the approximation spaces $V_h$ and $ Q_h$ are defined by
\begin{eqnarray}
V_h &=& \left\{ \bold{v}_h \in C^0(\Omega)^2 | \left. \bold{v}_h \right|_{K_h,i} \in \mathbb{P}_1^2, \forall i=1,\cdots,N-1  \right\},  \nonumber \\
Q_h &=& \left\{ q_h \in C^0(\Omega) | \left. q_h \right|_{K_{2h,j}} \in \mathbb{P}_1, \forall j=1,...,M-1 \right\}. \nonumber
\end{eqnarray}
Then, the velocity and the pressure are written
\begin{eqnarray}
p_h(x)=\sum_{j=1}^M p_{j}\phi_{j}
, \quad
\bold{u}_h(x)=\sum_{i=1}^N
\left(\begin{array}{ccc}
u_i\\
w_i
\end{array}
\right)\varphi_{i}.
\end{eqnarray}

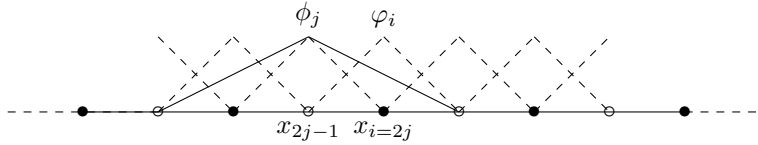
\begin{figure}[htbp]
\begin{center}
\begin{tikzpicture}
\draw[scale=1.0,domain=0:2,smooth,variable=\x] plot ({\x},{0.5*\x});

\draw[scale=1.0,domain=2:4,smooth,variable=\x] plot ({\x},{-0.5*\x+0.5*4});

\draw[dashed][scale=1.0,domain=0:1,smooth,variable=\x] plot ({\x},{\x});
\draw[dashed][scale=1.0,domain=1:2,smooth,variable=\x] plot ({\x},{\x-1});
\draw[dashed][scale=1.0,domain=2:3,smooth,variable=\x] plot ({\x},{\x-2});
\draw[dashed][scale=1.0,domain=3:4,smooth,variable=\x] plot ({\x},{\x-3});
\draw[dashed][scale=1.0,domain=4:5,smooth,variable=\x] plot ({\x},{\x-4});
\draw[dashed][scale=1.0,domain=5:6,smooth,variable=\x] plot ({\x},{\x-5});

\draw[dashed][scale=1.0,domain=0:1,smooth,variable=\x] plot ({\x},{-\x+1});
\draw[dashed][scale=1.0,domain=1:2,smooth,variable=\x] plot ({\x},{-\x+2});
\draw[dashed][scale=1.0,domain=2:3,smooth,variable=\x] plot ({\x},{-\x+3});
\draw[dashed][scale=1.0,domain=3:4,smooth,variable=\x] plot ({\x},{-\x+4});
\draw[dashed][scale=1.0,domain=4:5,smooth,variable=\x] plot ({\x},{-\x+5});
\draw[dashed][scale=1.0,domain=5:6,smooth,variable=\x] plot ({\x},{-\x+6});

    \draw [dashed] (-2,0) -- (0,0);
    \draw (-1,0) node[] {$\bullet$}
    	   -- (0,0) node[] {$\circ$}
       -- (1,0) node[] {$\bullet$}
       -- (2,0) node[] {$\circ$}
       -- (3,0) node[] {$\bullet$}
       -- (4,0) node[] {$\circ$}
       -- (5,0) node[] {$\bullet$}
       -- (6,0) node[] {$\circ$}
       -- (7,0) node[] {$\bullet$};
    \draw [dashed] (7,0) -- (8,0);
    \draw (2,1) node[above] {$\phi_{j}$} ;
    \draw (2,0) node[below] {$x_{2j-1}$} ;
    \draw (3,1) node[above] {$\varphi_{i}$} ;
    \draw (3,0) node[below] {$ x_{i=2j}$} ;
  \end{tikzpicture}
\end{center}
\caption{Representation of the basis functions.}
\label{figp1iso}
\end{figure}
In figure \ref{figp1iso}, the dashed lines are the usual elementary basis functions of $\mathbb{P}_1$ on the mesh $ K_h$, while the continuous lines are the basis functions on the mesh $ K_{2h}$. We can define the divergence operator, for all $j=1,M$
\begin{eqnarray}
\int_\Omega \div_{sw}(\bold{u}_h)\phi_j \, dx
= \sum_{K_h \in \mathcal{K}_h} \int_{K} \div_{sw}(\bold{u}_h)\phi_j \, dx.\nonumber
\end{eqnarray}
We use a linear interpolation for $H \varphi_j$, 
and consider that $ \Delta x_i = \Delta x \quad \forall i =1,\dots, N$ for the sake of simplicity. We
still approximate $\zeta$ by $\zeta_h$ defined before.

The discrete shallow water divergence operator is computed for all nodes $x_j$ of the mesh $\mathcal{K}_{2h} $ and therefore, denoting $ i=2j-1$, it can be written, $\forall j=1,M$
\begin{eqnarray*}
&& \int_\Omega \div_{sw}(\bold{u}_h)\phi_{j} \, dx
=
\left(\frac{1}{4} H_{i+2}{u}_{i+2}+\frac{H_{i+1}}{2}{u}_{i+1} \right)
-\left(\frac{1}{4} H_{i-2}{u}_{i-2}+\frac{H_{i-1}}{2}{u}_{i-1} \right) \nonumber \\
&&-\left(
\chi_{i-1/2} {m}_{i,j}^{{i-1/2}}
+
\chi_{i+1/2} {m}_{i,j}^{{i+1/2}}\right){u}_{i}
-\left(
\chi_{i-3/2} {m}_{i-2,j}^{{i-3/2}} \right){u}_{i-2}
-\left(
\chi_{i+3/2} {m}_{i+2,j}^{{i+3/2}} \right){u}_{i+2}
\nonumber\\
&&-\left(
\chi_{i-1/2} {m}_{i-1,j}^{{i-1/2}}
+
\chi_{i-3/2} {m}_{i-1,j}^{{i-3/2}}\right){u}_{i-1}
-\left(
\chi_{i+1/2}  {m}_{i+1,j}^{{i+1/2}}
+
\chi_{i+3/2}   {m}_{i+1,j}^{{i+3/2}} \right){u}_{i+1}
\nonumber\\
&&+ (2{m}_{i,j}^{{i+1/2}}){w}_i
+\left( {m}_{i-1,j}^{{i-1/2}}+{m}_{i-1,j}^{{i-3/2}}\right){w}_{i-1}
+\left( {m}_{i+1,j}^{{i+1/2}}+ {m}_{i+1,j}^{{i+3/2}}\right){w}_{i+1} \\
&&+({m}_{i-2,j}^{i-3/2}){w}_{i-2}
+{m}_{i+2,j}^{i+3/2}{w}_{i+2},
\end{eqnarray*}
with $ m_{i,j}^{i+1/2}=\int_{K_{h,i+1/2}}\varphi_{i}\phi_jdx$.

Similarly, the gradient shallow water operator is obtained for all the nodes $x_i$ of the mesh $ \mathcal{K}_h$. However,  we distinguish the gradient at the nodes of the elements $K_{2h}$ from the ones at the interior.
\noindent In other words, for all the nodes $ x_i$ of the mesh $\mathcal{K}_{2h}$, the gradient operator is defined  by
\begin{eqnarray*}
\left.\int_\Omega \nablasw p_h\cdot \underline{\varphi}_{(i=2j-1)} \, dx \right|_1&=&
\frac{H_{i}}{4}\left( p_{j+1}-p_{j-1}\right) \nonumber \\
&&+\chi_{i-1/2}
\left({m}_{i,j}^{i-1/2}p_j
+{m}_{i,j-1}^{i-1/2}p_{j-1}\right) \nonumber\\
&&+\chi_{j+1/2}\left(
{m}_{i,j}^{i+1/2}p_j
+{m}_{i,j+1}^{i+1/2}p_{j+1}\right),\\
\left.\int_\Omega \nablasw p_h\cdot\underline{\varphi}_{(i=2j-1)} \, dx \right|_2&=&
-2{m}_{i,j}^{i+1/2}p_j
 - {m}_{i,j}^{i-1/2}p_{j-1}
 -{m}_{i,j}^{i+1/2}p_{j+1} .
\end{eqnarray*}
On the other hand, for all the nodes $x_i$ such that $i$ is even
 \begin{eqnarray*}
\left.\int_\Omega \nablasw p_h\cdot\underline{\varphi}_{(i=2j)} \,dx \right|_1&=&
 \frac{H_{i}}{2}\left( p_{j+1}-p_{j}\right) \nonumber \\
&&
+\left(\chi_{i-1/2} {m}_{i,j}^{i-1/2}
      +\chi_{i+1/2} {m}_{i,j}^{i+1/2}
\right)p_{j} \nonumber\\
&&
+\left(\chi_{i-1/2} {m}_{i,j+1}^{i-1/2}
      +\chi_{i+1/2} {m}_{i,j+1}^{i+1/2}
\right)p_{j+1} ,\\
\left.\int_\Omega \nablasw p_h\cdot\underline{\varphi}_{(i=2j)} \, dx  \right|_2&=&
-\left( {m}_{i,j}^{i-1/2}+{m}_{i,j}^{i+1/2}\right)p_{j}
-\left( {m}_{i,j+1}^{i-1/2}+{m}_{i,j+1}^{i+1/2}\right)p_{j+1}.
\end{eqnarray*}

With the discretization of the shallow water operators given below, we are able to validate the scheme for the two first order methods. Nevertheless, notice that in the following section, only the first order has been implemented.

\section{Validation with an analytical solution}
\label{sec:valid}

In \cite{JSM_nhyd_num},\cite{JSM_nhyd} some analytical solutions of the model \eqref{nh1}-\eqref{contrainte} have been presented and they allow to validate the numerical method. We consider the propagation of a solitary wave without topography. This solution has the form
\begin{eqnarray*}
H &=& H_0+a\left(\hbox{sech} \left(\frac{x-c_0 t}{l}\right)\right)^2, \\
u &=& c_0 \left( 1-\frac{d}{H} \right),\\
w &=& -\frac{a c_0 d}{l H} \hbox{sech} \left( \frac{x-c_0 t}{l} \right) \hbox{sech}'\left( \frac{x-c_0 t}{l} \right), \\
p &=& \frac{a c_0^2 d^2}{2l^2 H^2}\left( \left( 2H_0 -H\right)\left( \hbox{sech}'\left(\frac{x-c_0 t}{l} \right) \right) ^2 \right. ,\\
& & + \left. H\hbox{sech}\left( \frac{x-c_0 t}{l}\right) \hbox{sech}''\left( \frac{x-c_0 t}{l}\right)\right),
\end{eqnarray*}
with  $d,a,H_0 \in \mathbb{R} $, $ H_0>0$, $a>0 $
and $ c_0=\frac{l}{d} \sqrt{\frac{gH^3}{l^2-H_0^2}}$, $l = \sqrt{\frac{H_0^3}{a}+H_0^2} $. \\

\begin{figure}[htbp]
\begin{center}
\begin{tabular}{cc}
\includegraphics[scale=0.3]{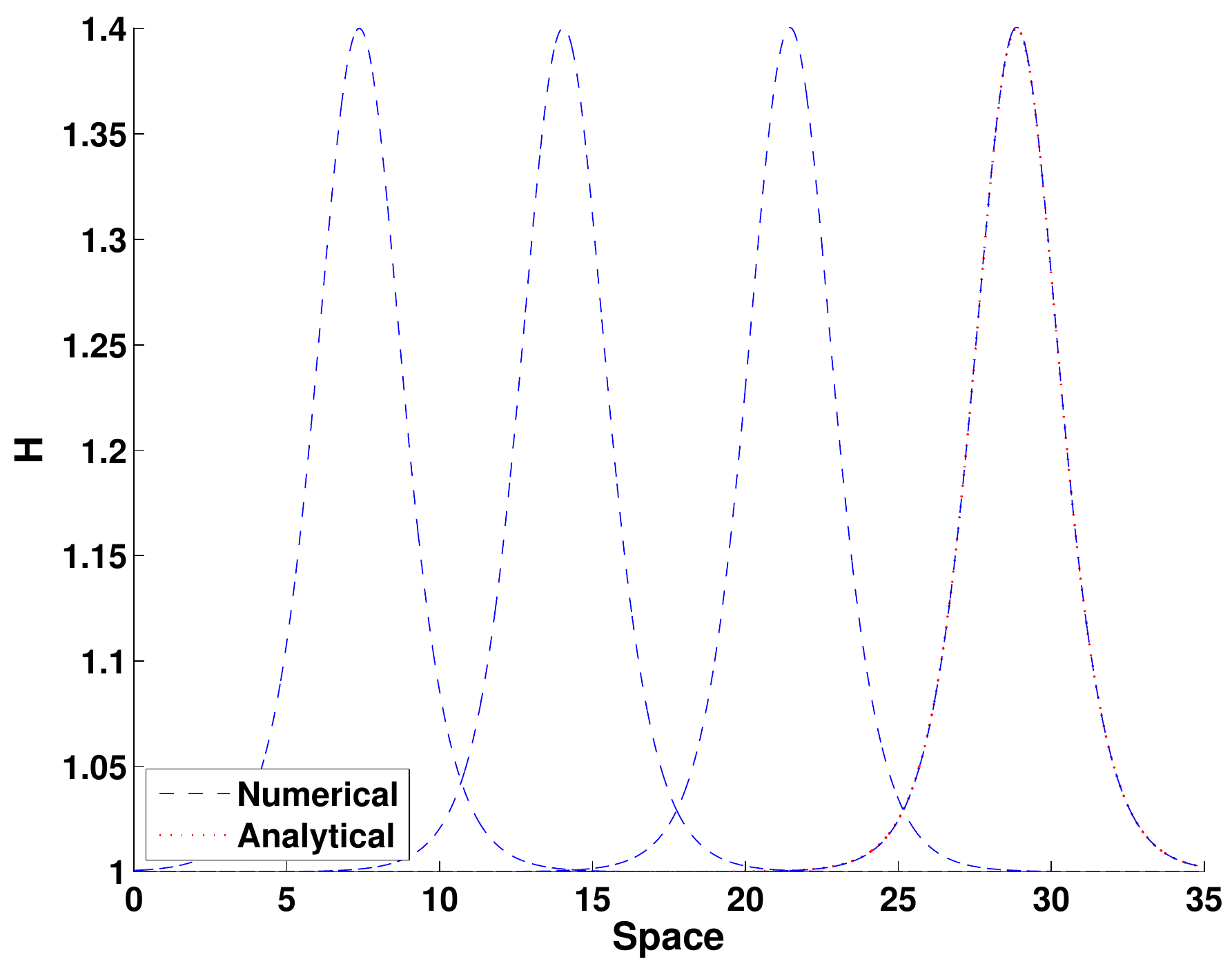} &
\includegraphics[scale=0.3]{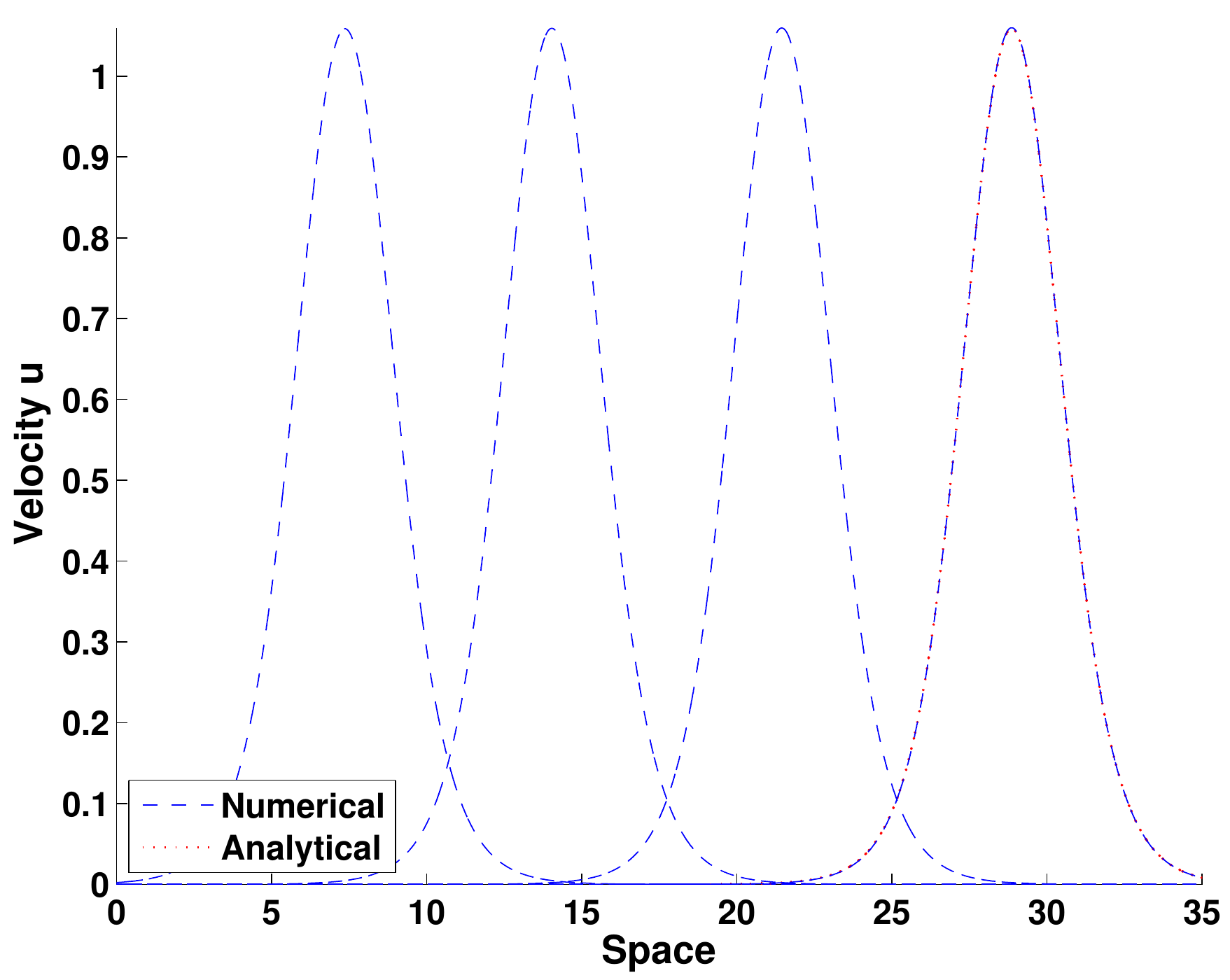}\\

\includegraphics[scale=0.3]{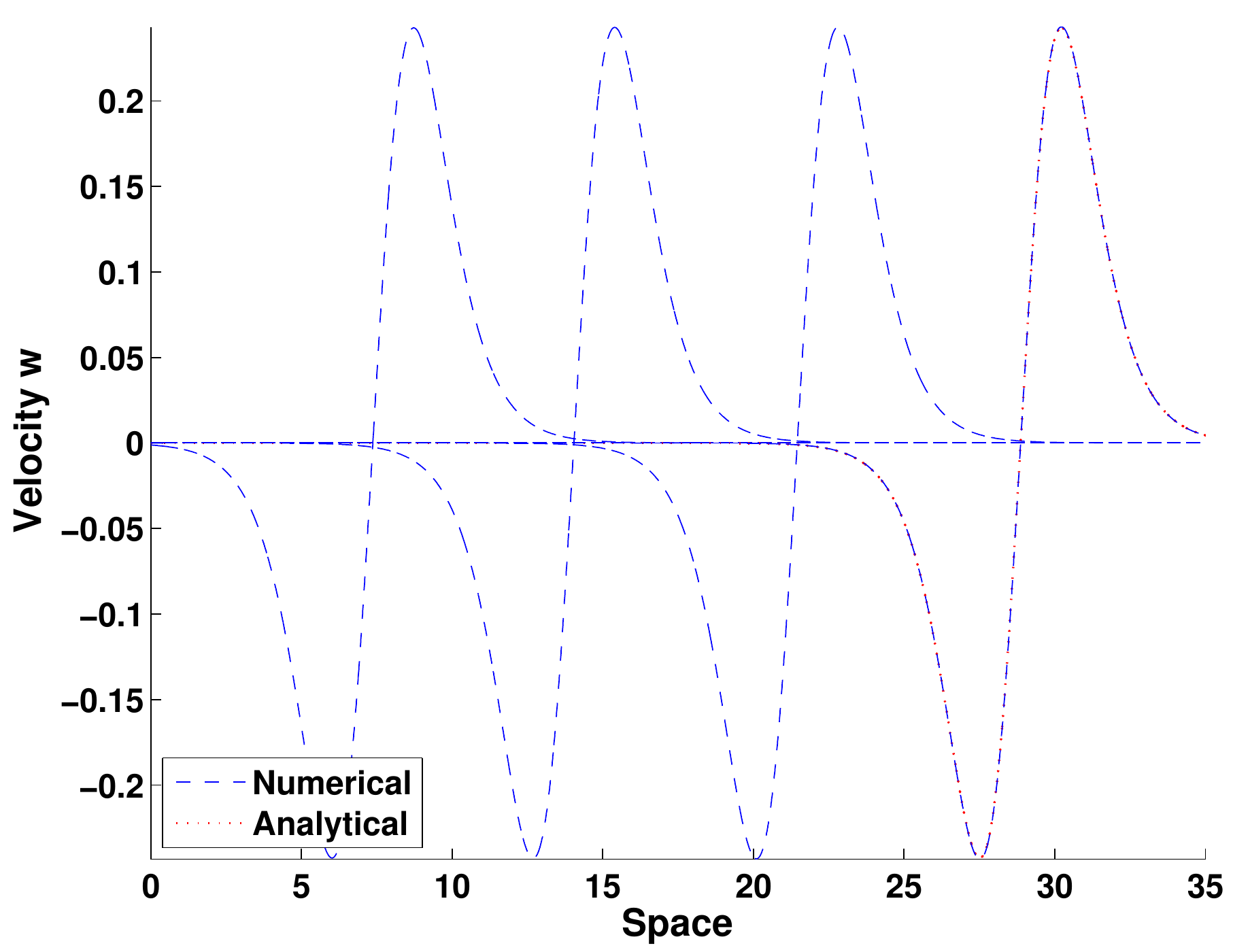} &
\includegraphics[scale=0.3]{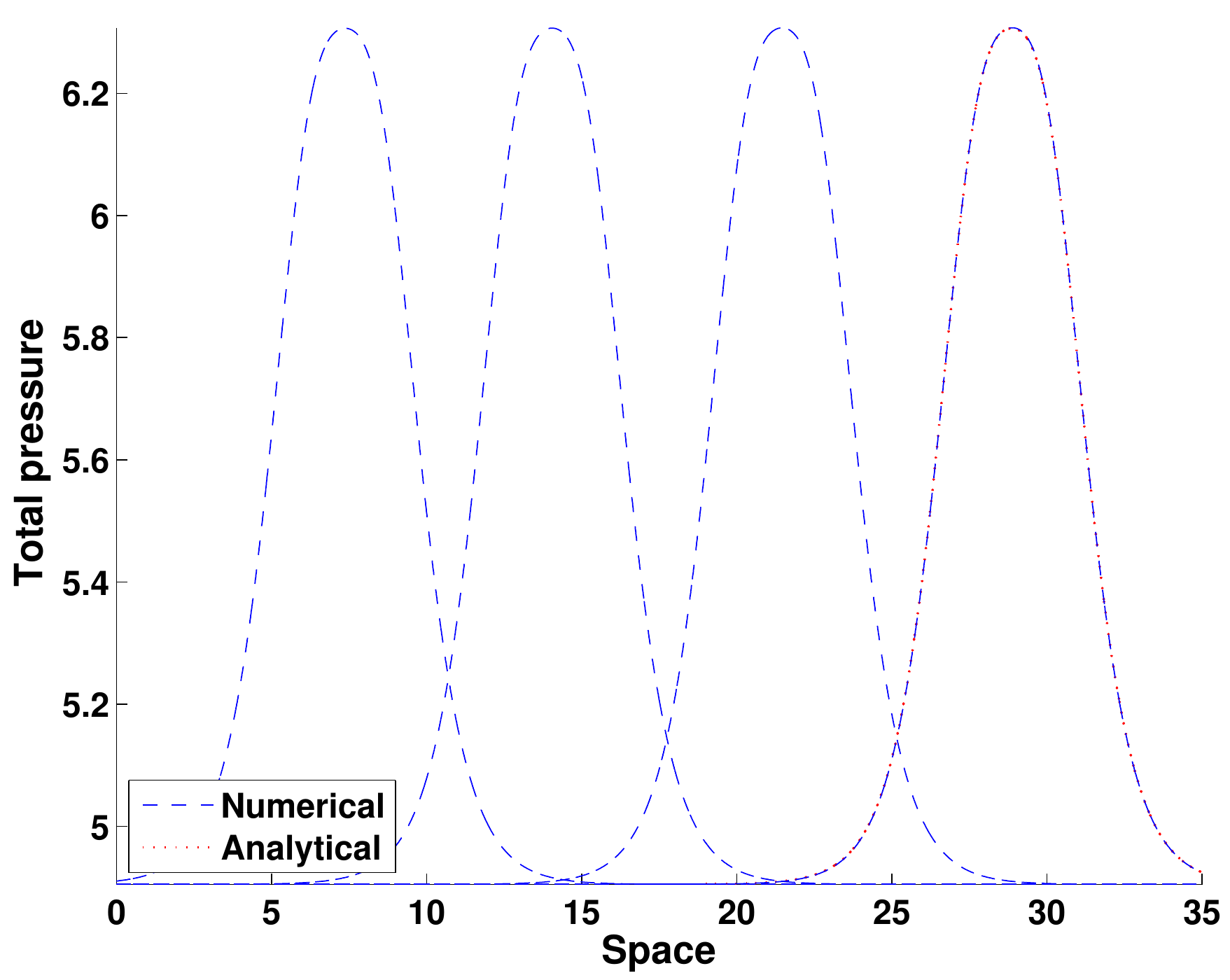}
\end{tabular}
\caption{Propagation of the solitary wave at times  $ 1.00008\, s ,
  1.9009\, s, 3.9017\, s$ and $5.9025\, s$. Comparison with analytical solution at time $t= 5.9025 \, s $.}
\label{figSoliton}
\end{center}
\end{figure}

The solitary wave is a particular case where dispersive contributions are counterbalanced by non linear effects so that the shape of the wave remains unchanged during the propagation. The propagation of the solitary wave has been simulated for the parameters $a=0.4\,m, H_0=1\,m$, and $ d=1\,m$ over a domain of $ 45 \,m$ with $9000$ nodes.
At time $t=0$, the solitary wave is positioned inside the domain. The results presented in figure \ref{figSoliton} show the different fields, namely the elevation, the components of velocity and the total pressure at different times, and the comparison with the analytical solution at the last time.

In the projection step, the greatest difficulty is to compute the pressure corresponding to the boundary conditions of the hyperbolic part (as seen in \ref{secBD}). The solution near the boundary has been confronted to the analytical solution. In the following result, we set a Neumann boundary condition on the non hydrostatic pressure with the parameters given below. As shown in the figure \ref{soliton_pnh}, the pressure is well estimated at the outflow boundary and allows the wave to leave the domain with a good behavior.
 \begin{figure}[htbp]
\begin{center}
\includegraphics[scale=0.4]{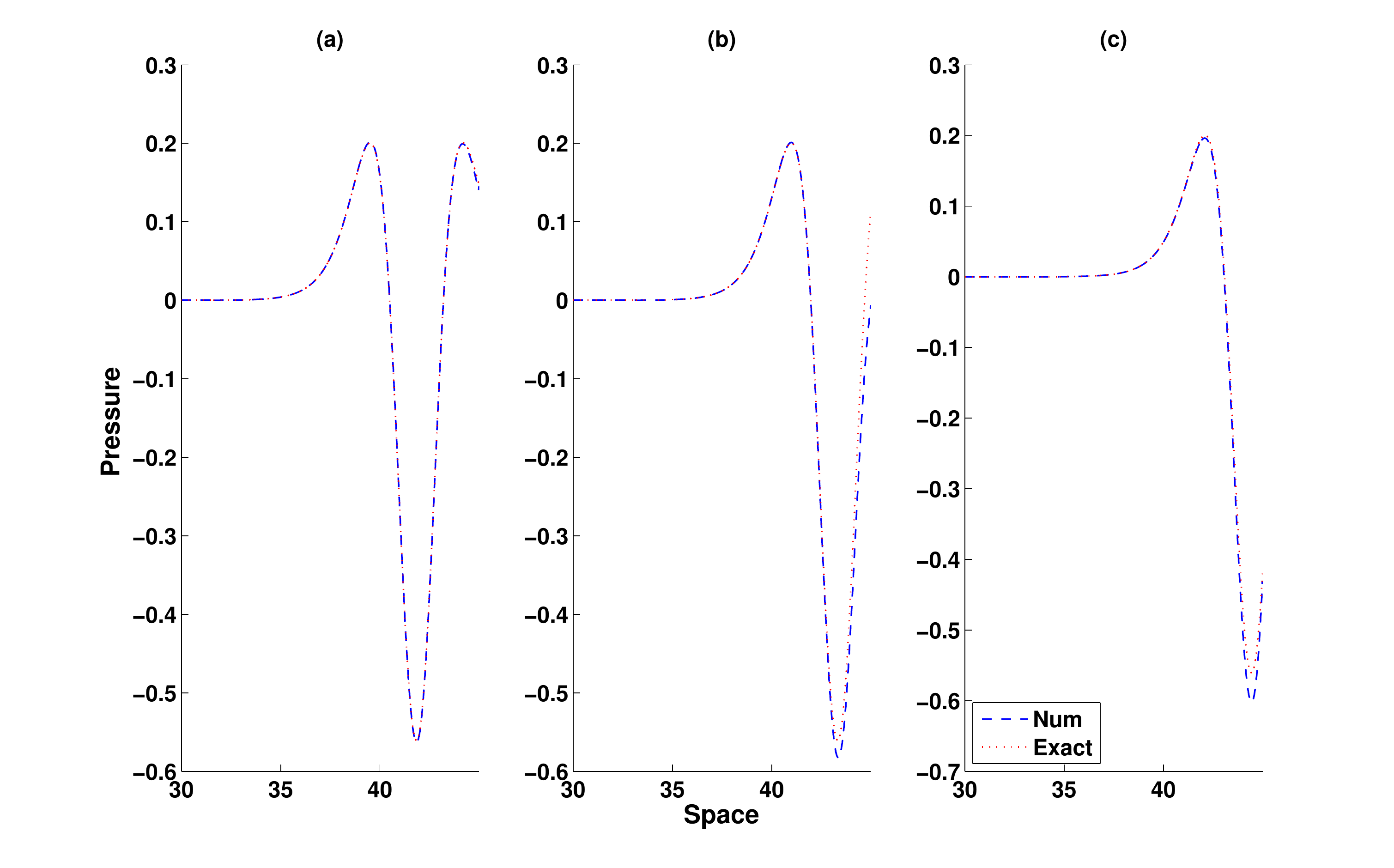}
\caption{Non hydrostatic pressure profile at right boundary $(x=45\, m)$
(a):
$t=9.4044 \,s$
(b):
$t=9.8046\, s$
(c):
$t=10.1048 \,s$.}
\end{center}
\label{soliton_pnh}
\end{figure}
The inflow boundary condition has been tested with this same test case and gives similar results. We are able to let the solitary wave enter in the domain with a good approximation of the elevation.

\begin{figure}[htbp]
\begin{center}
\includegraphics[scale=0.4]{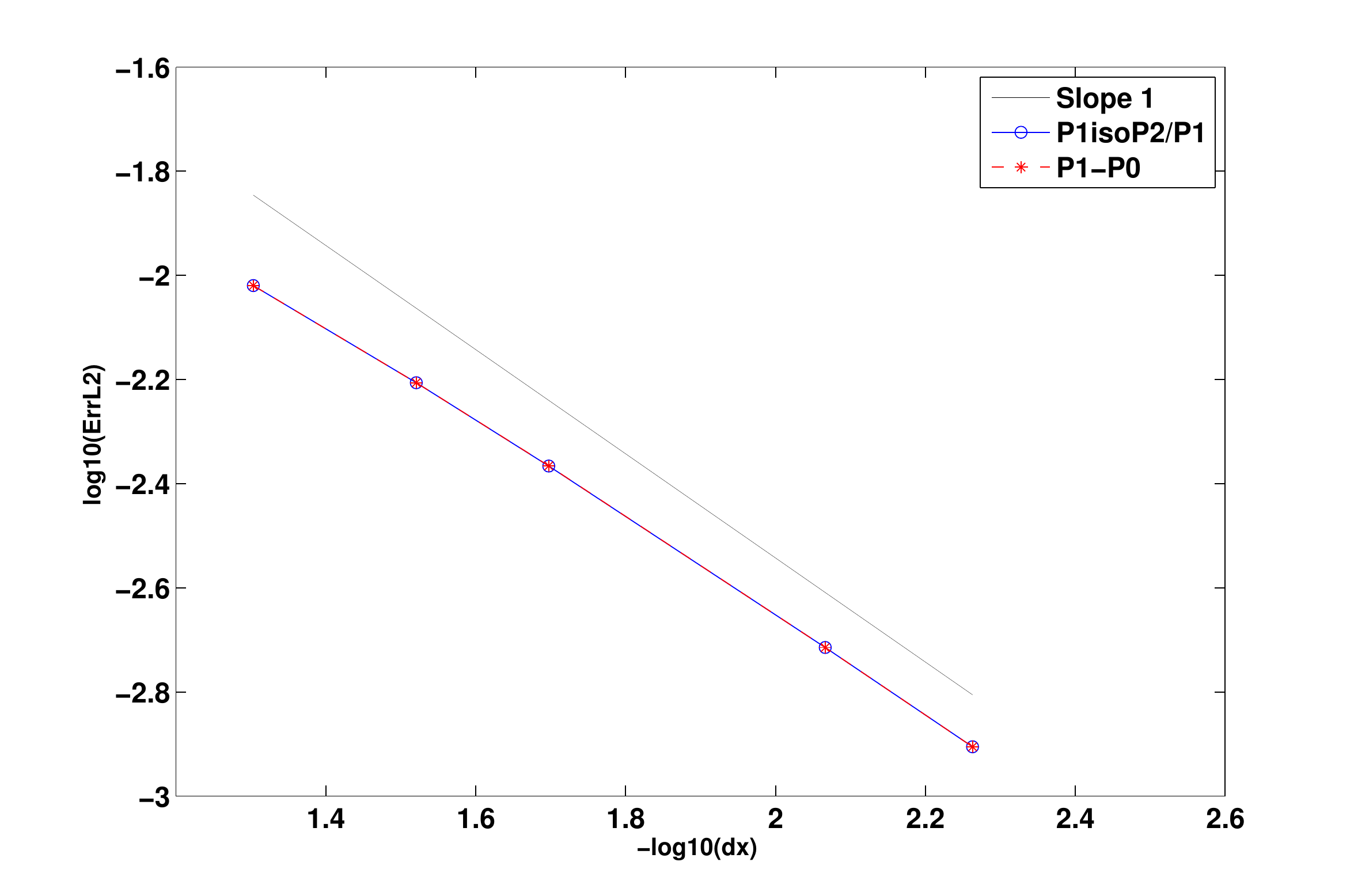}
\caption{Convergence rate of the solitary wave solution for $\mathbb{P}_1$/$\mathbb{P}_0$ and $\mathbb{P}_1$-iso-$\mathbb{P}_2$/$\mathbb{P}_1$  scheme.}
\label{figOrder}
\end{center}
\end{figure}

The numerical simulations for the first order method are compared with the analytical solution and the $ L^2 $- error has been evaluated over different meshes of sizes from 603 nodes to 6495 nodes (see figure \ref{figOrder}). With the parameters given above, it gives a convergence rate close to $1$ for the two computations, i.e $ \mathbb{P}_1$/$\mathbb{P}_0$ and $\mathbb{P}_1$-iso$\mathbb{P}_2$/$\mathbb{P}_1$.

Notice that the parameters set to validate the method lead to have a
significant non hydrostatic pressure (see the figure
\ref{soliton_pnh}) and then, the results show the ability of the
method to preserve the solitary wave over the time. The numerical results have also been obtained for the Thacker's test presented in \cite{JSM_nhyd_num}, with the same rate as the $\mathbb{P}_1$/$\mathbb{P}_0$ method.

\section{Numerical results}
\subsection{Dam break problem}

\begin{figure}[htbp]
\begin{center}
\begin{tabular}{cc}
\includegraphics[scale=0.2]{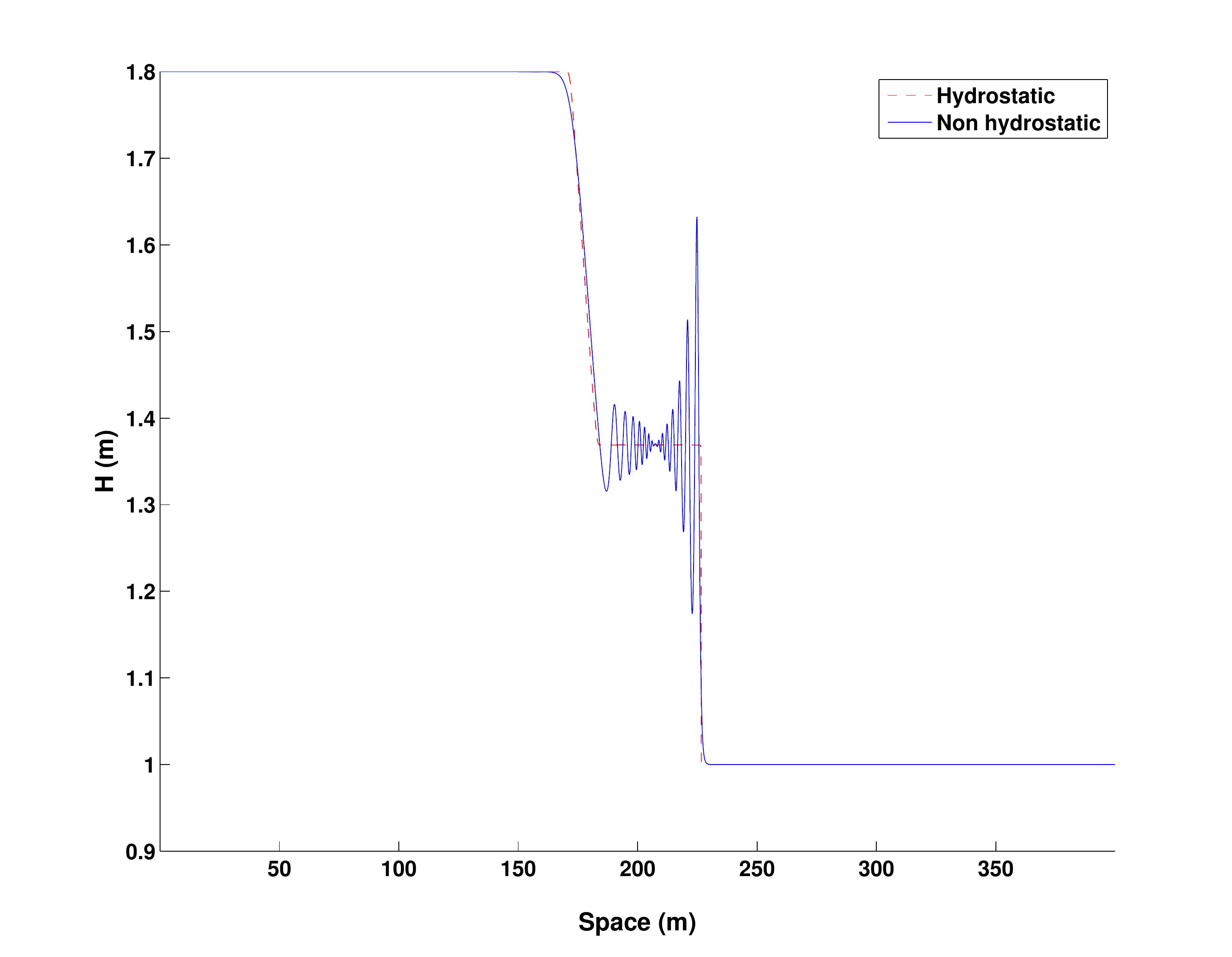} &
\includegraphics[scale=0.2]{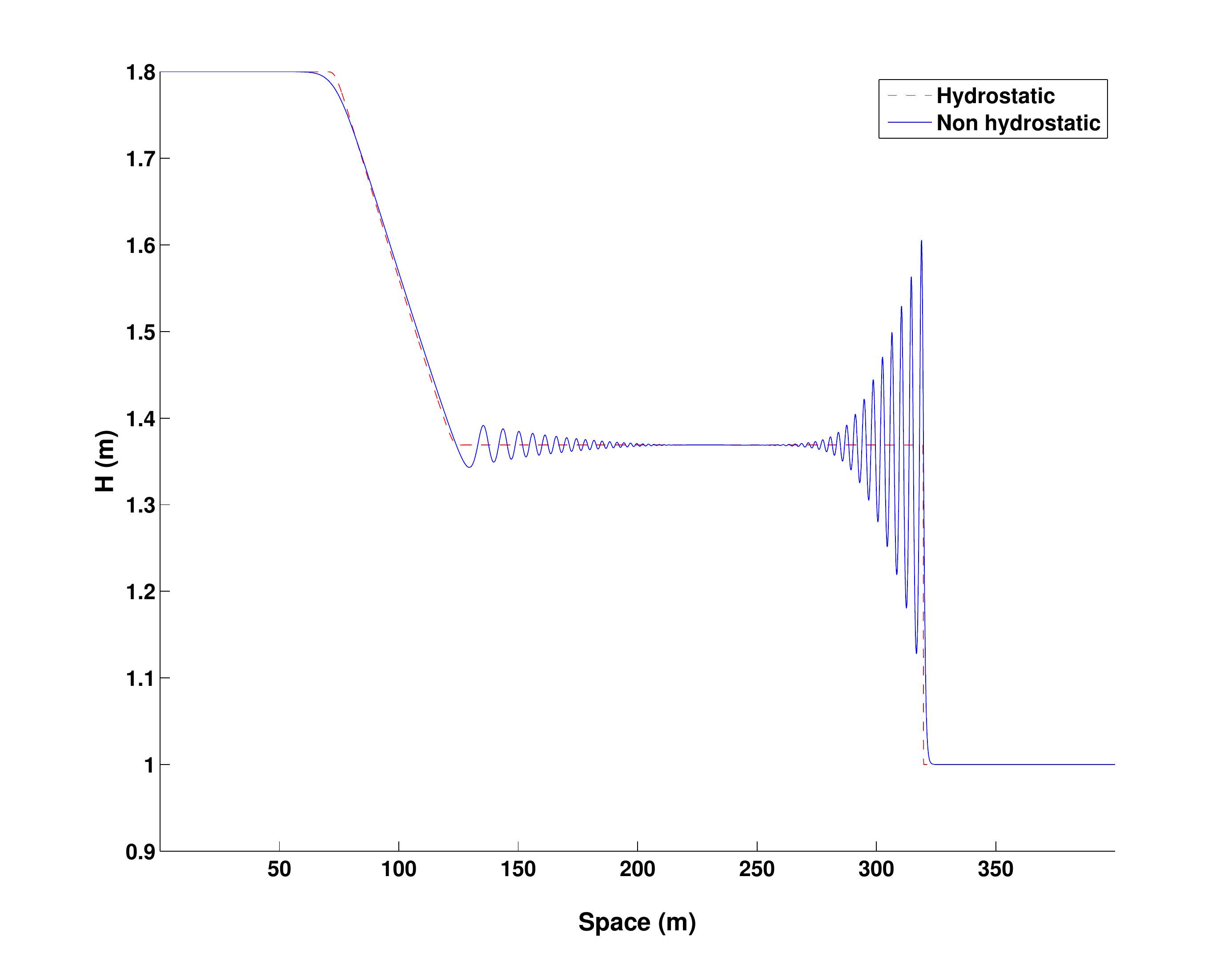}\\
\includegraphics[scale=0.2]{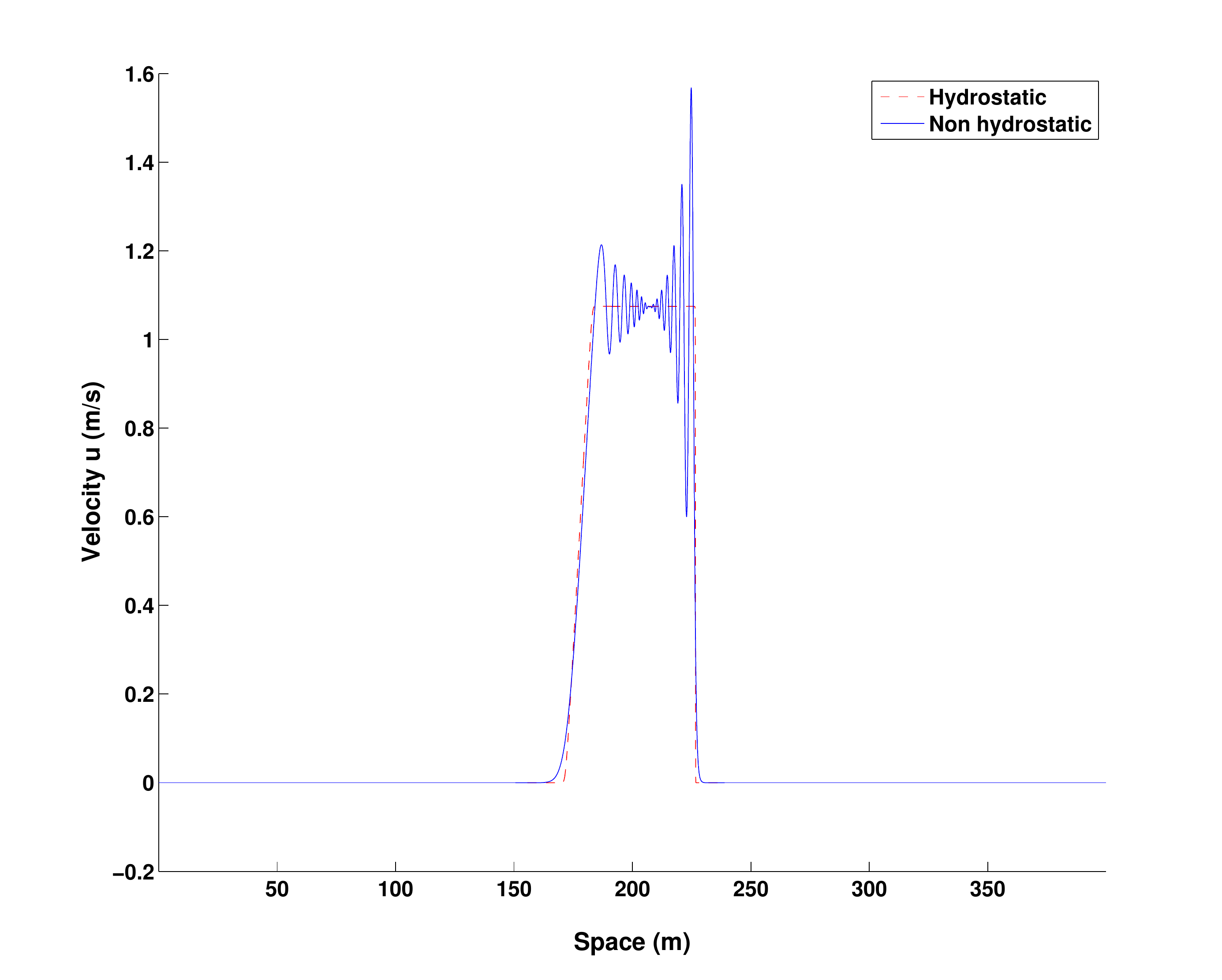} &
\includegraphics[scale=0.2]{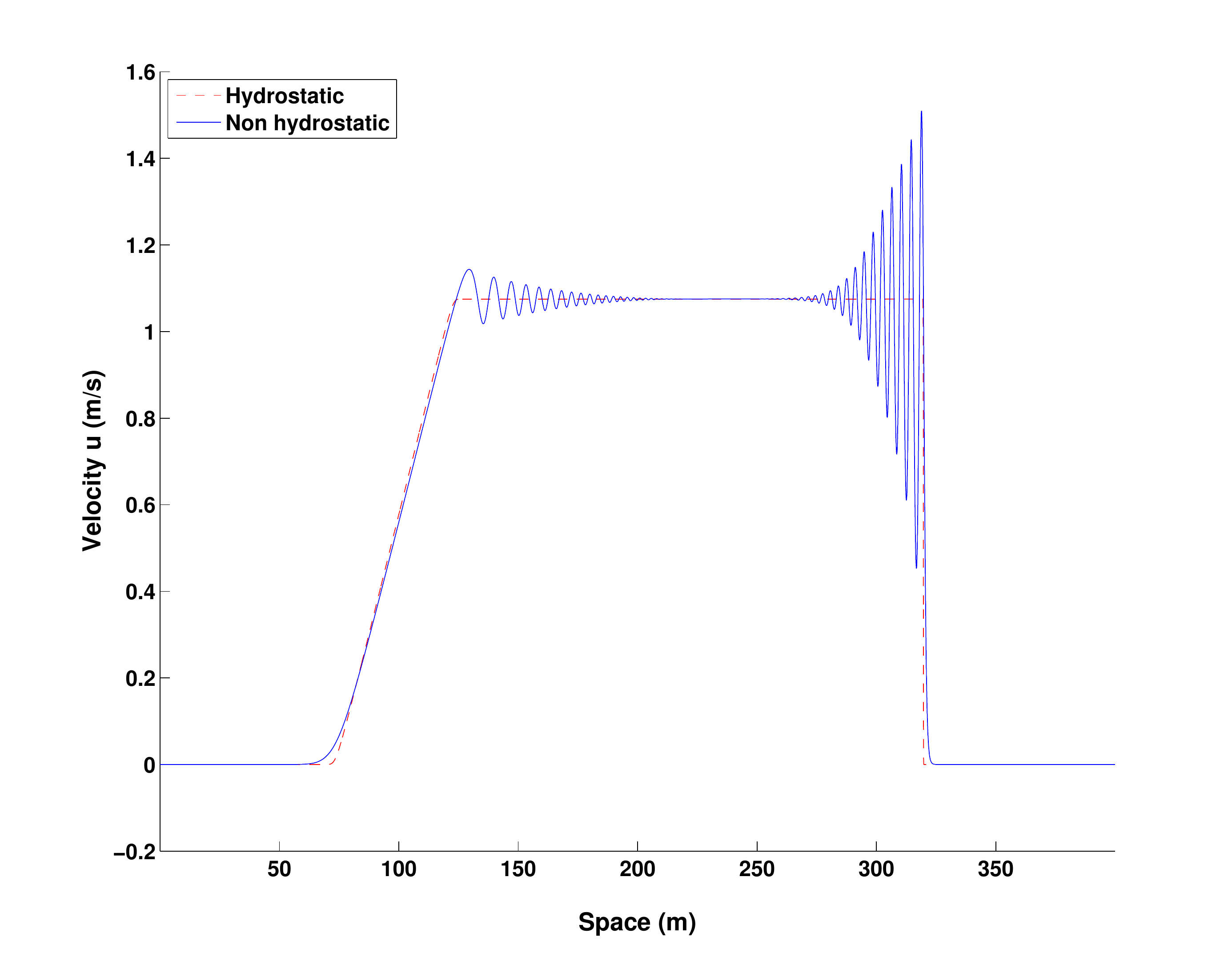}
\end{tabular}
\caption{The dam break problem, elevation $H$ and velocity $u$ at times $ t=10 \,s$ and  $ t= 45 \, s$.}
\label{figDam}
\end{center}
 \end{figure}

We next study the dispersive effect on the classical dam break
problem, which is usually modeled by a Riemann problem providing a
left state $(H_L,u_L)$ and a right state $(H_R,u_R)$ on each side of
the discontinuity $x_d$ (\cite{godlewski_book}). However, our
numerical dispersive model does not allow discontinuous solutions due to the variational spaces required for $H$ (see also \cite{JSM_nhyd}), thus we provide an initial data numerically close to the analytical one
\begin{eqnarray*}
H(x,0) &=& (H_R+a)-a\tanh\left(\frac{x-x_d}{\epsilon}\right),\\
a &=& H_R-H_L.
\end{eqnarray*}
To evaluate the non hydrostatic effect, the different fields have been compared with the shallow water solution with the initial data:
$H_L = 1.8\, m$, $H_R = 1\, m$, $u_R = u_L = 0\, m.s^{-1}$,
 $\epsilon = 10^{-4}\, m$, $x_d = 300\, m$ over a domain of length $600 \, m $ with $30 000$ nodes.
%
%
In figure \ref{figDam}, the evolution of the state is shown at time
$t= 10\, s$ and $t = 45 \,s $. The oscillations are due to the
dispersive effects but the mean velocity does not change. These results are in adequation with the analysis proposed by Gavrilyuk in \cite{lemet-gavri} for the Green-Naghdi model with the same configuration.

\subsection{Wet-dry interfaces}

The ability to treat the wet/dry interfaces is crucial in geophysical problems, since geophysicists are interested in studying the behavior of the water-depth near the shorelines.
This implies a water depth tending to zero at such boundaries. To treat the problem, we use the method introduced in \cite{JSM_nhyd_num}, considering a minimum elevation $H_\epsilon$.

\begin{center}
\begin{figure}[htbp]
\begin{center}
\includegraphics[scale=0.25]{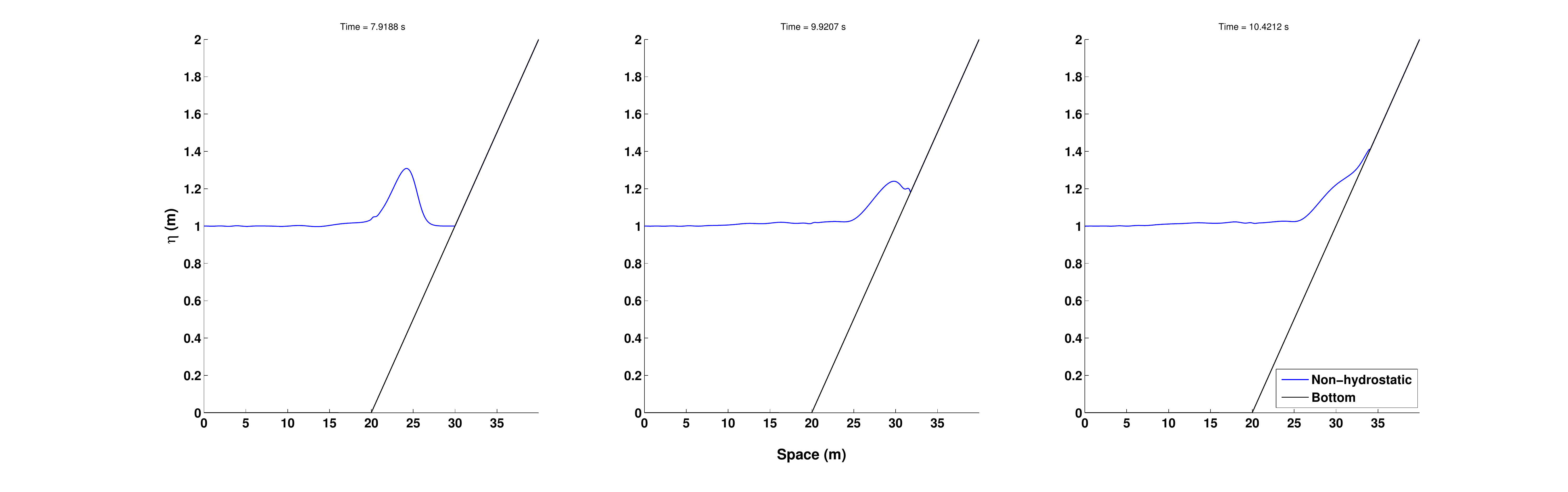}
\end{center}
\caption{Propagation of a wave at a wet/dry interface.}
\label{figWet}
\end{figure}
\end{center}

Therefore, we confront the method with a coastal bottom at the right
boundary over a domain of $ 35\, m$ with $3000 $ nodes. A wave is
generated at the left boundary with an amplitude of $ 0.2 \,m$ and an
initial water depth $ H_0=1\,m$. In figure~\ref{figWet}, the arrival
of the wave at the coast is shown for times $t=7.91\,s $, $ 9.92\,s$ and $ 10.42\,s $.

\subsection{Comparison with experimental results}

\begin{figure}[htbp]
\begin{center}
\includegraphics[scale=0.2]{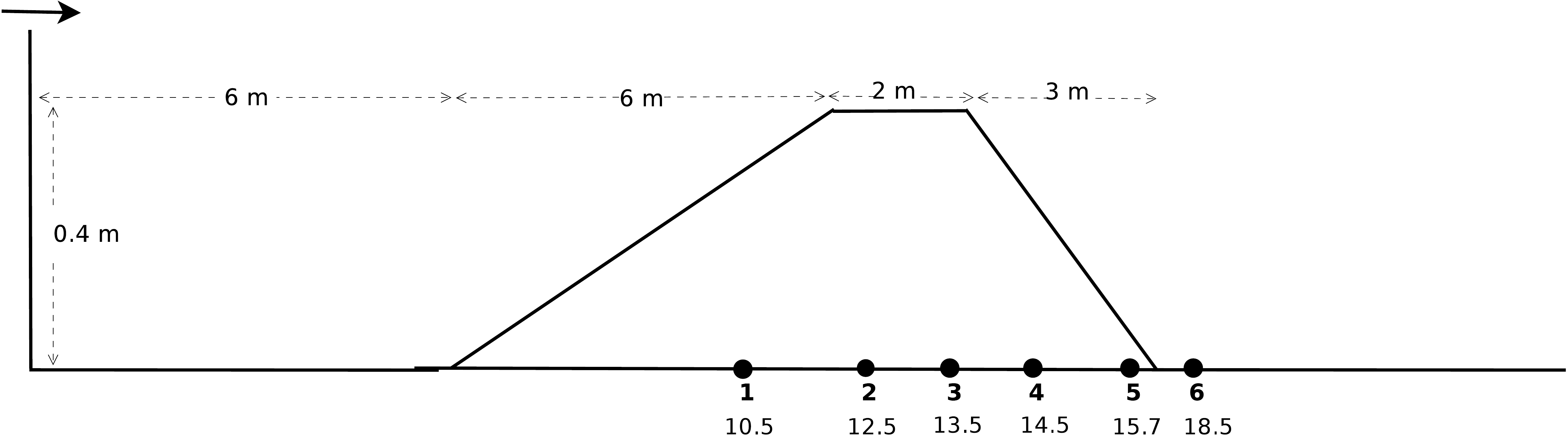}
\caption{Configuration of Dingemans's test.}
\label{geoDing}
\end{center}
\end{figure}

In this part, we confront the model with Dingemans experiments (detailed in \cite{DingemansTEst,dingemans}) that consist in generating a small amplitude wave at the left boundary of a channel with topography as described in figure~\ref{geoDing}.

\begin{figure}[htbp]
\begin{center}
\includegraphics[scale=0.35]{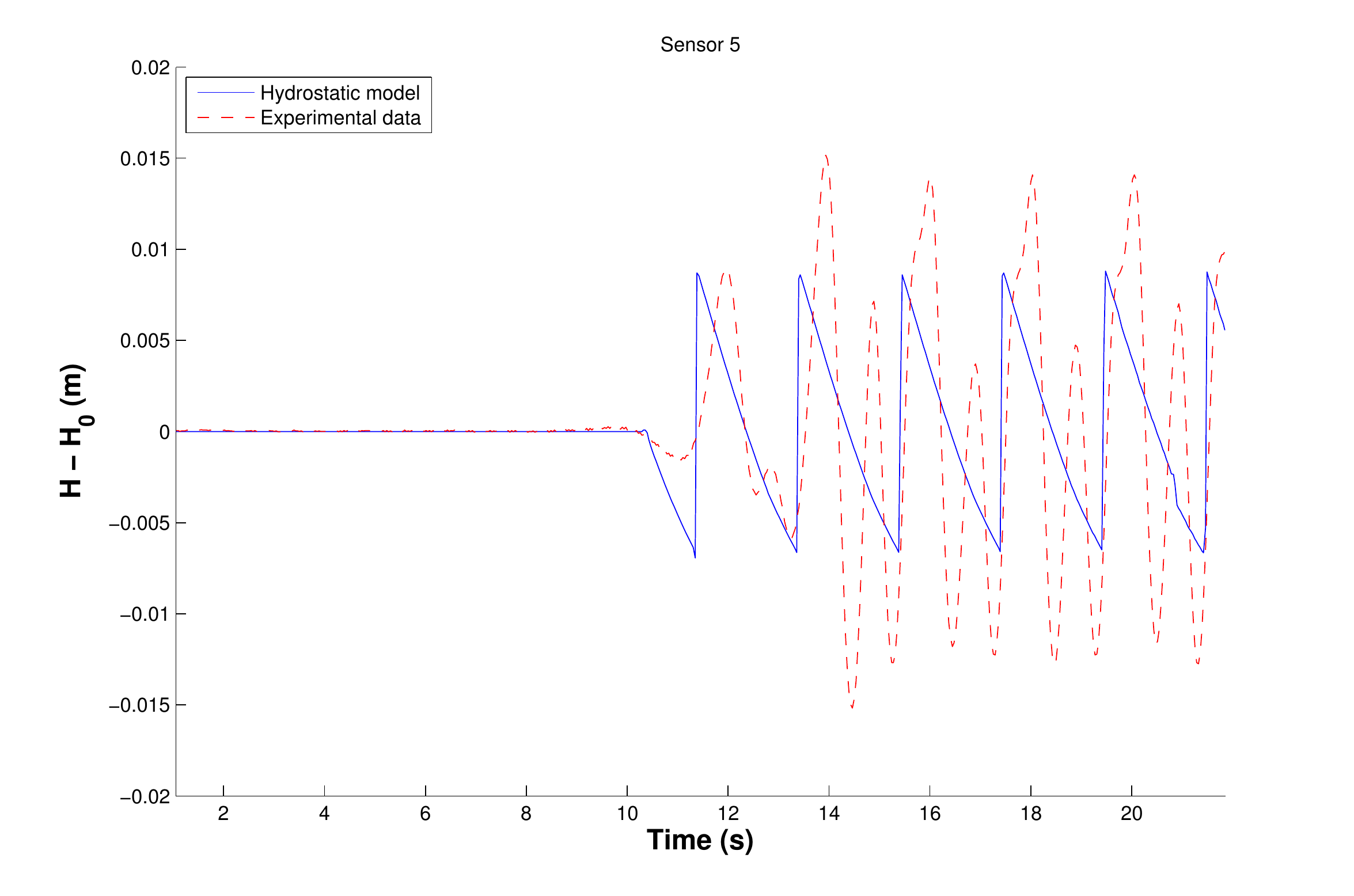}
\caption{Comparison with hydrostatic model on sensor 5.}
\label{dingHyd}
\end{center}
\end{figure}

At the left boundary, a wave is generated with a period $ T=2.02 \, s
$ and an amplitude of $ 0.02 \, m $. A free outflow condition is set
at the right boundary. The initial free surface is set to be
$\eta_0=0.4 \, m$, and the measurement readings are saved  at the
following positions $10.5\,m$, $12.5\,m$, $13.5\,m$, $14.5\,m$
,$15,7\,m$ and $17.3\,m$, placed at sensors $1$ to $6$ (fig.\ref{geoDing} ).
In such a situation, the non hydrostatic effects have a significant impact on the water depth that cannot be represented by a hydrostatic model. These effects result mainly from the slope of the bathymetry, $10\%$ in this case.
In the figure \ref{dingHyd}, the simulation has been run with the hydrostatic model and the elevation has been compared with measures at the sensor $5$. As one can see, the non-hydrostatic pressure has to be taken in consideration to estimate the real water depth variation.

\begin{figure}[htbp]
\begin{center}
\includegraphics[scale=0.3]{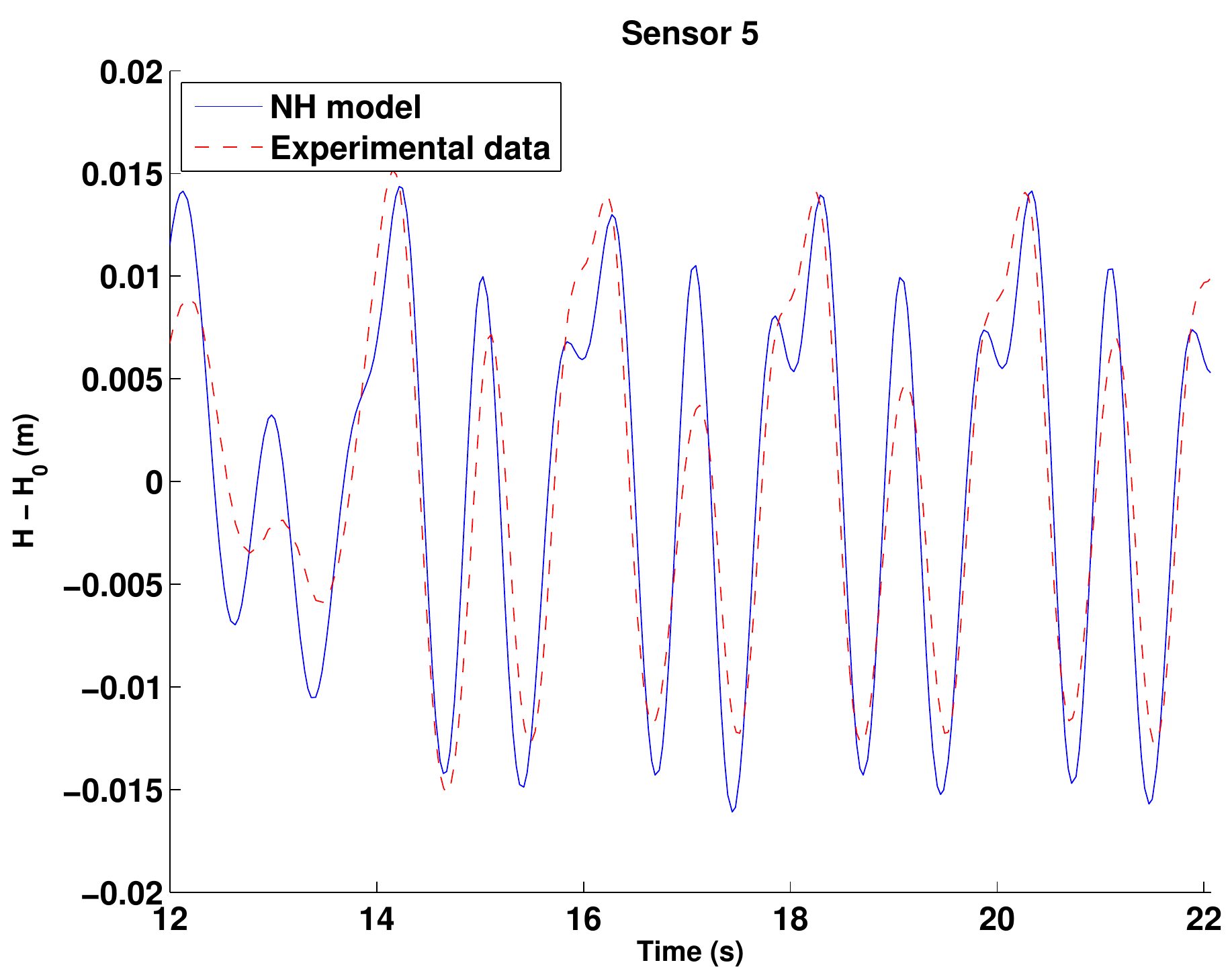}
\includegraphics[scale=0.3]{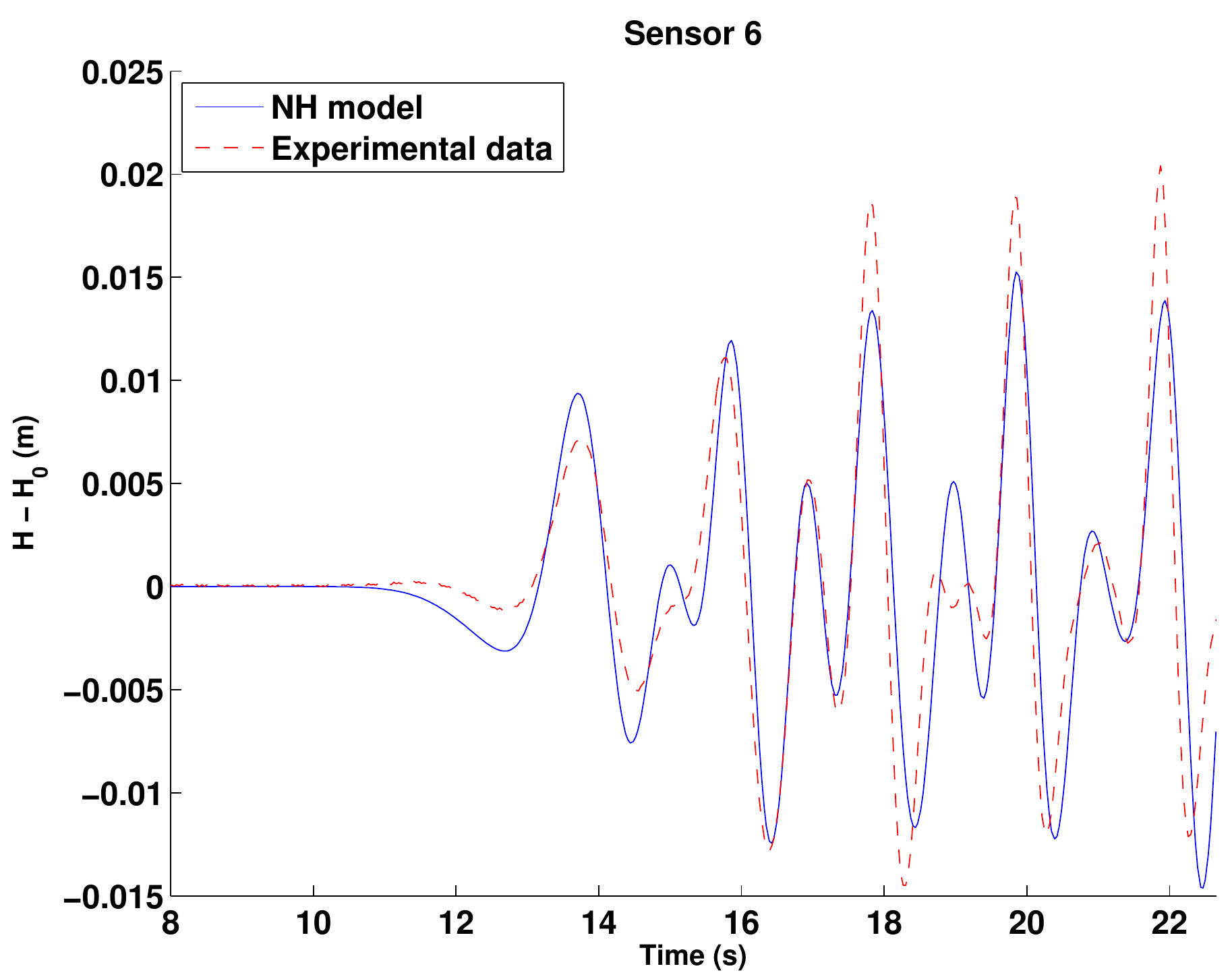}
%
\caption{Comparison between measured and computed elevations on
  Dingemans test for the first six sensors.}
\label{figDing}
\end{center}
\end{figure}

The numerical simulation with the non-hydrostatic model has been run with $ 15 000 $ nodes on a domain of $ 49\, m$ over $ 25\,s$ and the comparisons are illustrated for each sensor (fig. \ref{figDing}).
\begin{figure}[htbp]
\begin{center}
\includegraphics[scale=0.3]{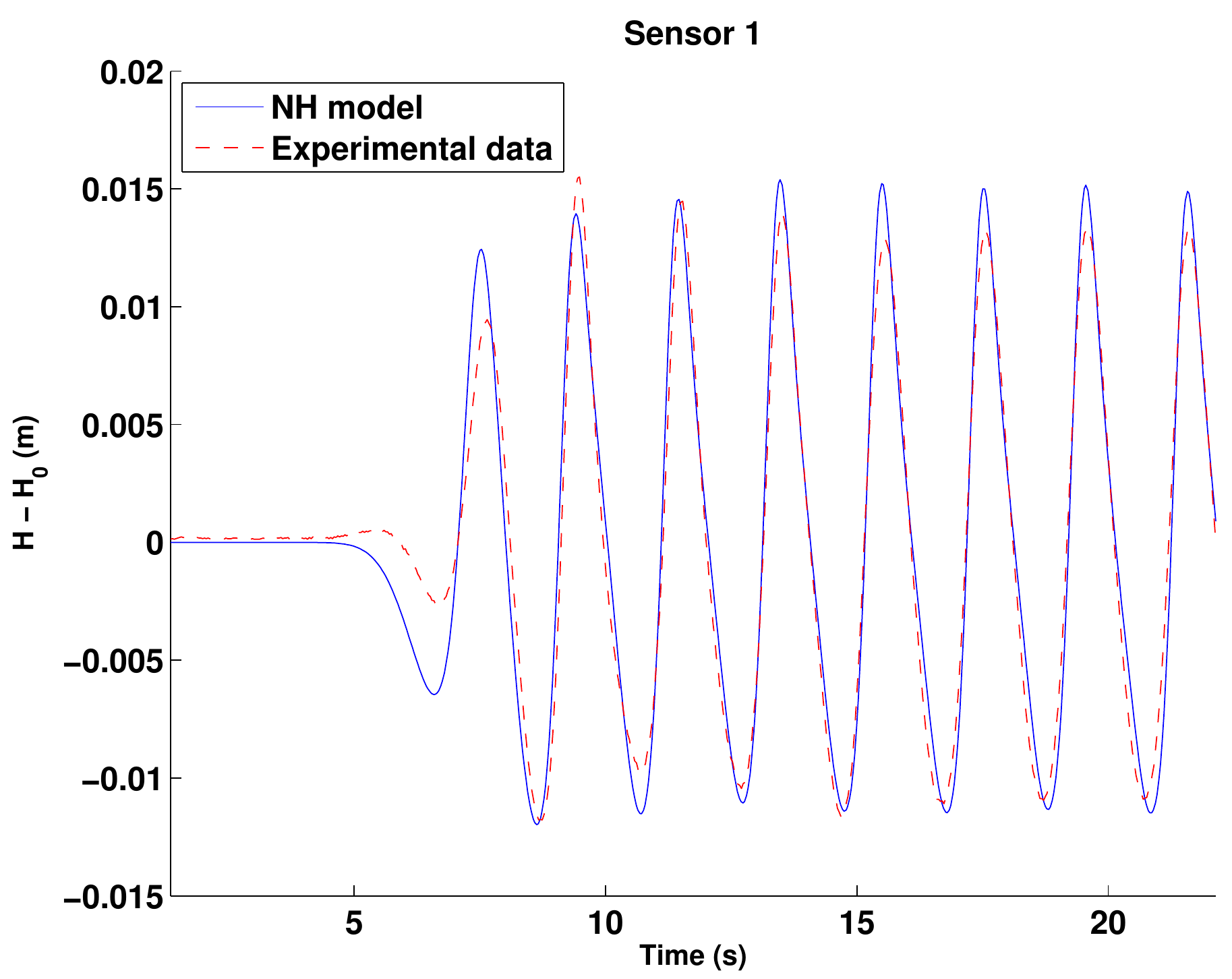}
\includegraphics[scale=0.3]{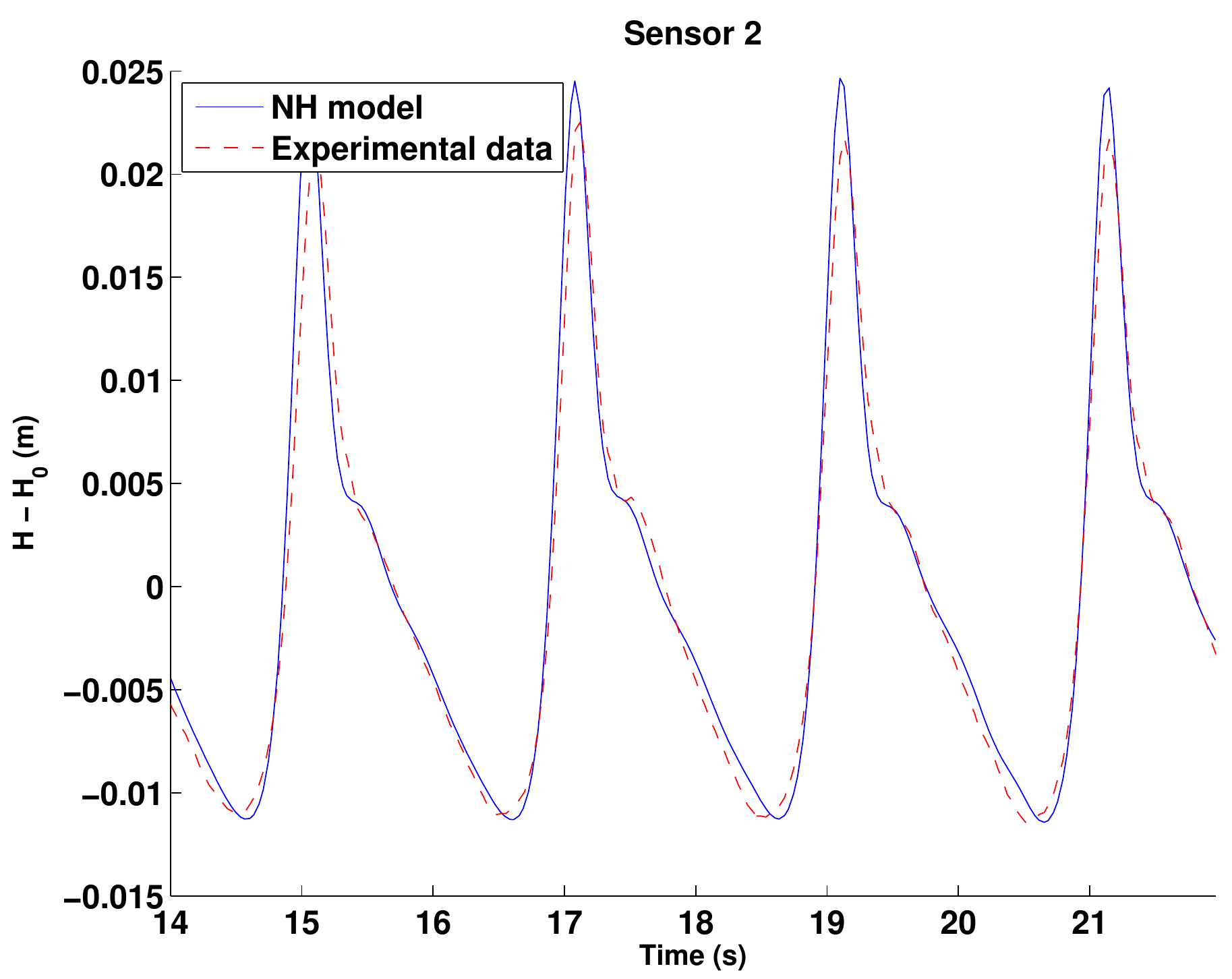}
\end{center}
\end{figure}
\begin{figure}[htbp]
\begin{center}
\includegraphics[scale=0.3]{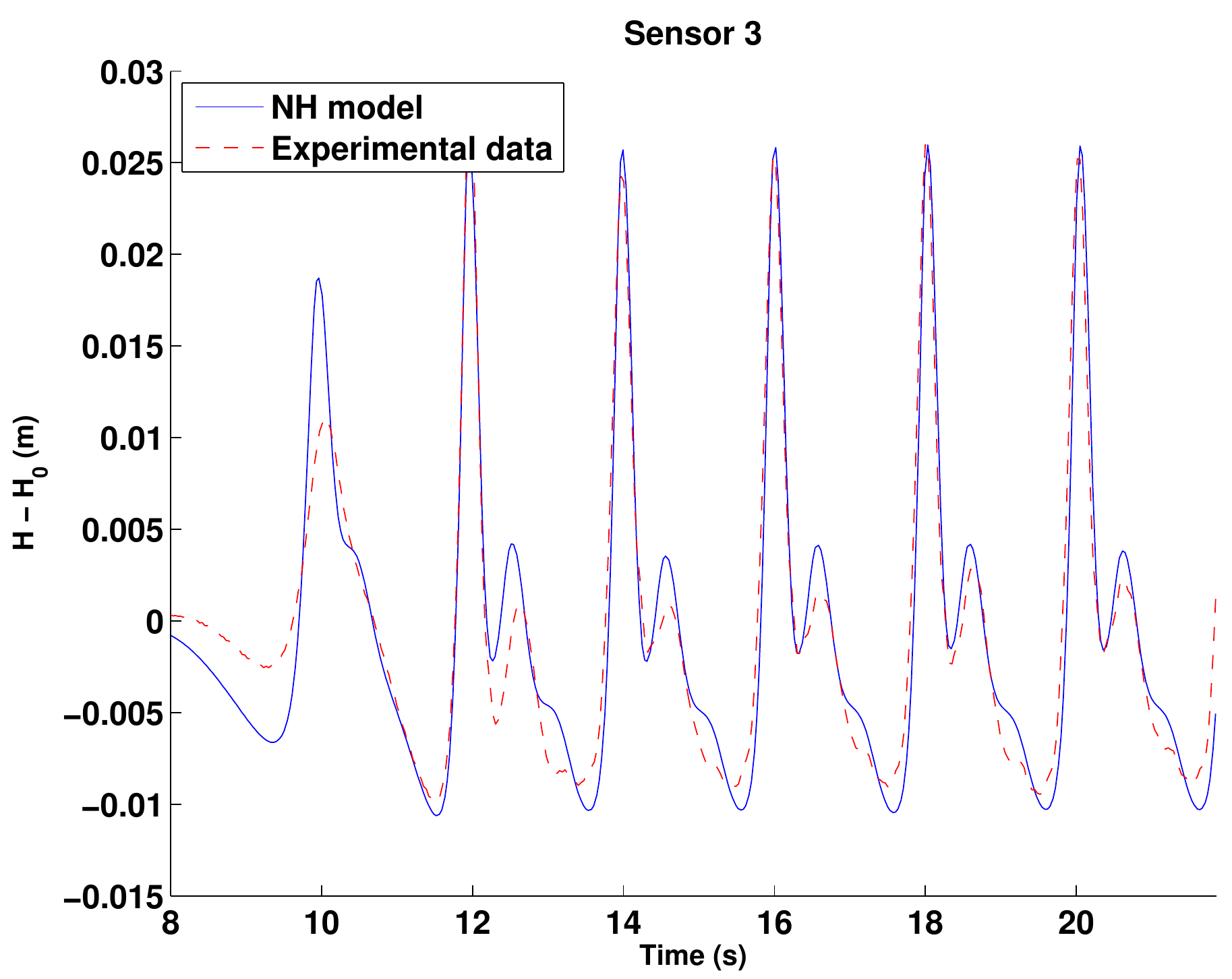}
\includegraphics[scale=0.3]{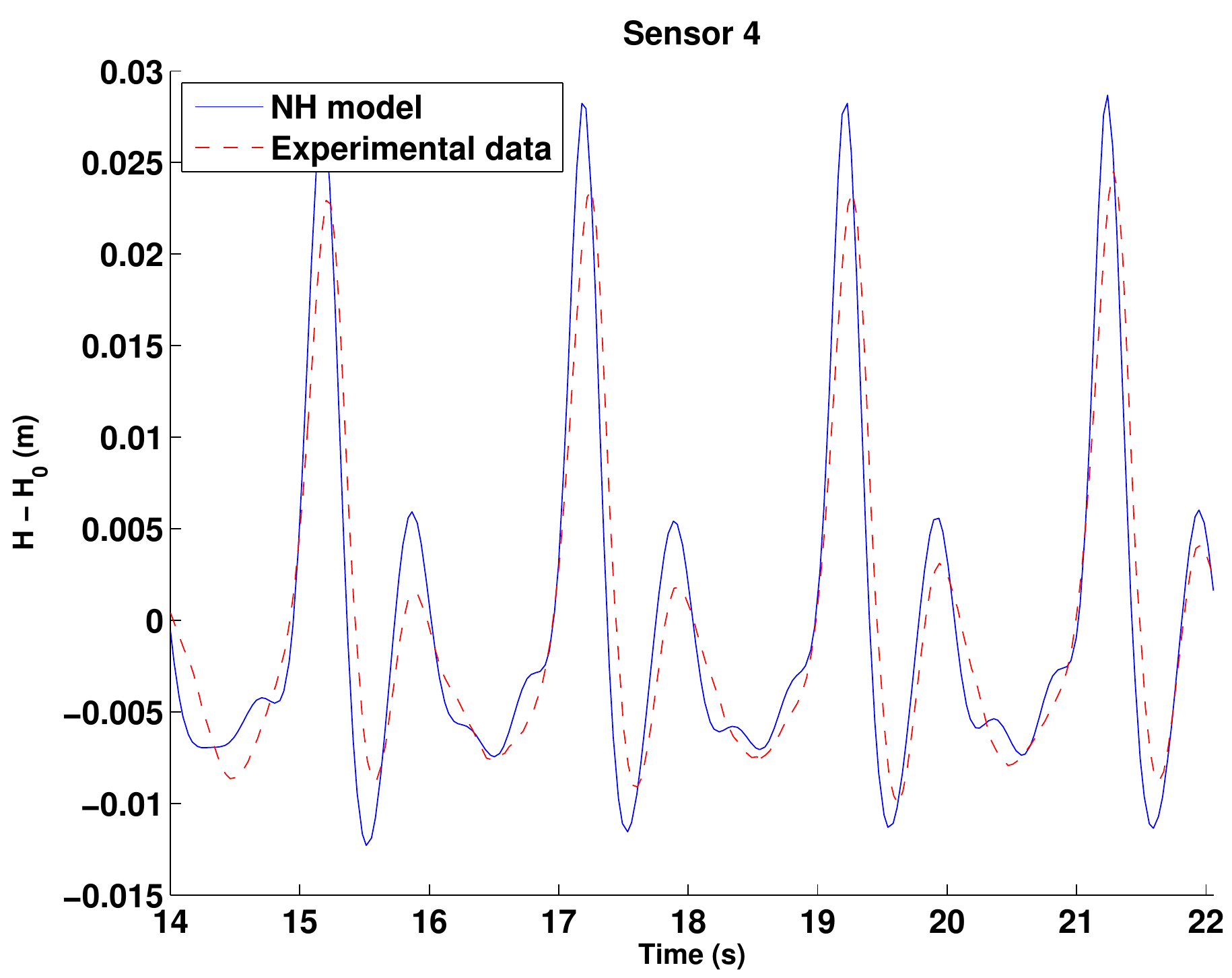}
\end{center}
\end{figure}
The goal of this last result is also to highlight the ability of the model to capture dispersive effects for a geophysical flow with a non negligible pressure.

\subsection{Remark on iterative method}

\begin{figure}[htbp]
\begin{center}

\includegraphics[scale=0.45]{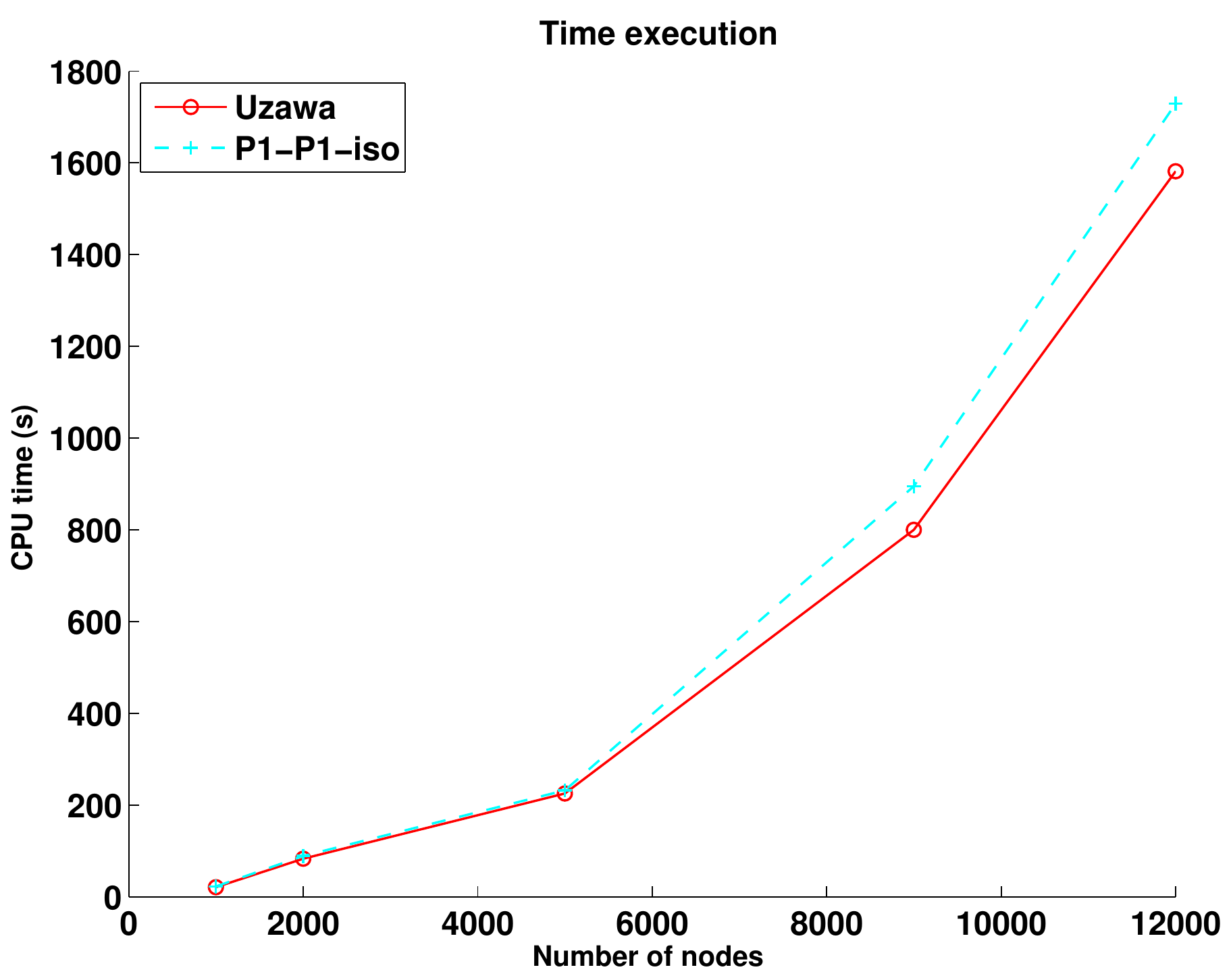}
\caption{Comparison of the computing time (CPU) for the $ \mathbb{P}_1$/$\mathbb{P}_0$ scheme, the $ \mathbb{P}_1$-iso-$\mathbb{P}_2/\mathbb{P}_1$ and Uzawa method with $ \mathbb{P}_1$-iso-$\mathbb{P}_2/\mathbb{P}_1$ solver with a tolerance $10^{-5}$.}
\label{figtime}
\end{center}
\end{figure}

We recall that this
formulation should allow to extend the method on two dimensional
unstructured grids. However, it requires to inverse a system at each
time iteration, which will become too costly in two dimensions. To anticipate the two dimensional problem, this method has been tested using different iterative methods like conjugate gradient and Uzawa methods. In figure \ref{figtime}, we show a comparison of the computing time for the implementation $\mathbb{P}_1$-iso-$\mathbb{P}_2/\mathbb{P}_1$ and Uzawa method. In one dimension, it is not relevant to use one of these methods, while it will be necessary for the two dimension model.

\section{Conclusion}

In this paper, a variational formulation has been established for the one dimensional dispersive model introduced in \cite{JSM_nhyd}. The main idea is to give a new framework in which it will be possible to extend the scheme to the two dimensional model. To this aim, the finite-element method has been presented with two approximation spaces. First, the $\mathbb{P}_1$/$\mathbb{P}_0$ approximation has been done and we recover, as expected, the finite difference scheme, together with the good results proved in \cite{JSM_nhyd_num}. Then, the  $\mathbb{P}_1$-iso-$\mathbb{P}_2$/$\mathbb{P}_1$approximation has been studied to prepare the two dimensional problem.
We have validated the method using several numerical tests and studying the dispersive effect on geophysical situations.

\bibliographystyle{plain}
\bibliography{biblioA}

\end{document}